\documentclass[letter,12pt]{article}
\usepackage{amsmath}
\usepackage{enumerate}
\usepackage{amssymb}
\usepackage{geometry}
\usepackage{setspace}
\usepackage{graphicx}
\usepackage{natbib}
\usepackage{lmodern}
\usepackage{soul}
\usepackage{comment}
\usepackage{ifthen}
\usepackage{xr}
\usepackage{framed}
\newtheorem{theorem}{Theorem}

\newtheorem{definition}{Definition}
\newtheorem{lemma}{Lemma}
\newtheorem{proposition}{Proposition}
\newtheorem{remark}{Remark}

\newenvironment{proof}[1][Proof]{\textbf{#1.} }{\ \rule{0.5em}{0.5em}}
\DeclareMathSizes{12}{11}{7.7}{5.5}
\usepackage[inline]{enumitem}
\usepackage{comment}

\DeclareMathOperator*{\esssup}{ess\,sup}
\DeclareMathOperator*{\argmin}{argmin} 

\graphicspath{ {WD-Graphs/} }

\begin{document}


\thispagestyle{empty}
\title{Solutions to Equilibrium HJB Equations for Time-Inconsistent Deterministic Linear Quadratic Control: Characterization and Uniqueness }
\author{ Yunfei Peng \and Wei Wei 
\thanks{ Peng: School of Mathematics and Statistics, Guizhou University, Guiyang 550025, China. Email: yfpeng@gzu.edu.cn. Wei: Department of Actuarial Mathematics and Statistics, Heriot-Watt University, Edinburgh, Scotland EH14 4AS, UK. E-mail: wei.wei@hw.ac.uk.  } }
%
%
\thispagestyle{empty}
\date{\today}
\maketitle
\thispagestyle{empty}

\vspace*{-0.4cm}
\begin{abstract}This paper investigates a class of HJB equations that delineate equilibria for time-inconsistent deterministic linear-quadratic control problems within the intra-personal game theoretic framework. The source of inconsistency is non-exponential discount functions. Our approach involves characterizing solutions to the HJB equations through a specific class of Riccati equations, incorporating integral terms. Through a thorough examination of the uniqueness of solutions to these integro-differential Riccati equations, we establish the uniqueness of solutions to the equilibrium HJB equations.
\end{abstract}


\newpage
\setcounter{page}{1}
\section{Introduction}\label{SecIntroduction}

Time-inconsistency describes a situation where optimal policies determined in the past are no longer considered optimal today and, therefore, are not implemented. Individual decision-makers with time-inconsistent preferences are often referred to as ``dividuals'', signifying the existence of different selves at different points in time who may not act consistently among themselves. Consequently, dynamic decision making becomes time-inconsistent, hence time-inconsistent control problems arise.
As is standard in the literature on decision making, time-inconsistent control problems are often considered within the intra-personal game theoretic framework and the corresponding equilibria are taken as solutions to such problems.\footnote{See, for example, \cite{phepol1968}; \cite{lai1997}; \cite{odorab1999}; \cite{krusmi2003} and \cite{lutmar2003}.}\par

Behavioral scientists and economists have been studying time-inconsistency problems for more than sixty years\footnote{In economics, \cite{str1955} has observed that non-constant time preference rates result in time-inconsistent decisions.}. However, serious mathematical treatments to time-inconsistent control problems had not been available until around 15 years ago. Particularly, \cite{bjomur2009} introduce a systematic method to derive HJB equations for general time-inconsistent control problems. Since then, HJB equations have been widely used to construct intra-personal equilibria for time-inconsistent models in various areas\footnote{See, for example, \cite{harlai2013}, \cite{bjomurzho2014} and \cite{ebert2020weighted}.}. While most of literature on time-inconsistent control is focused on the existence of the solutions to the HJB equations,  the uniqueness is relatively unexplored. As far as we know, The paper of \cite{ekepir2008} is the only one that mentions a uniqueness result for the HJB equations, in which a non-uniqueness result is reported in a time-inconsistent portfolio management problem.\par

In this paper we establish the uniqueness of solutions to the HJB equations for a general deterministic time-inconsistent linear-quadratic (LQ) control problem, where the inconsistency arises from non-exponential discount functions\footnote{Non-exponential discount functions are a major source of time-inconsistency. This type of discount functions would stem from decreasing impatience, time-insensitivity and collective decision making, etc. See \cite{tha1981}, \cite{ebepre2007} and \cite{weitzman2001gamma} for more details.}. While classical, time-consistent LQ problems can be solved without resorting to dynamic programming, the time-inconsistent variants inherently require a dynamic, intrapersonal game-theoretic approach. In such contexts, the equilibrium HJB equation plays a central and indispensable role in characterizing desirable strategies. The uniqueness is challenging in time-inconsistent environment, as the comparison principle, which is almost the most powerful tool to establish the uniqueness of solutions to HJB equations, is rarely obtained in the inter-personal game theoretic framework. Moreover, the HJB equations arising in time inconsistent control problems usually link to the notion of feedback control (see Definition \ref{DefEquilibrium}). Due to loss of (global) optimality in time-inconsistent control problems, the involvement of feedback control would damage the linear structure of the equilibrium when we implement the spike variation on the LQ control problem\footnote{This feature marks a significant difference between equilibria based on notions of feedback control and open loop control in time-inconsistent problems. }. As a consequence, it would be difficult to know the structure of the value function by the quadratic form of the objective functional. This is distinct from the conventional time-consistent LQ control theory and constitutes the major difficulty in the uniqueness analysis of corresponding time-inconsistent HJB equations.

In order to overcome the difficulty, we develop a novel method to characterize the solutions to the time-inconsistent HJB equations without using the optimality  of the value function or the quadratic structure of the objective functional directly (See Proposition \ref{pl2}.). Our method is motivated by the equivalence between Riccati equations and control problems established by \cite{cai2022time}. The key observation is that an integral related to the second order derivative of the value function solves a time-inconsistent Riccati equation. This observation implies the quadratic structure of the value function and makes a linkage between the Riccati equation and the HJB equation. We then prove the uniqueness of the solution to the Riccati equation. Thanks to the uniqueness result for the Riccati equation, we finally obtain the uniqueness of the solution to the HJB equation for the LQ control problem.\par

Recent years have witnessed a rapid growth of the studies on time-inconsistent LQ control within the intra-personal game theoretic framework\footnote{\cite{he2022time} offers comprehensive overviews of the literature on time inconsistency in general control theory.}. Early research on this topic include \cite{basak2010dynamic}, \cite{bjomurzho2014}, \cite{hu2012time} and \cite{yong2012deterministic}. \cite{basak2010dynamic} and \cite{bjomurzho2014} study dynamic mean-variance portfolio selection problems which can be formulated within the time-inconsistent LQ framework. The equilibria in their papers are based on the spike variation and defined in the class of feedback policies. This definition of equilibria is formally proposed by \cite{ekelaz2006} and \cite{bjomur2009} and consistent with the intra-personal equilibria pursued in the majority of the literature on behavioral economics, in which the equilibria are outcomes of recursions\footnote{See \cite{bjomur2009}  for more discussions and examples.}. \cite{yong2012deterministic}, \cite{yong2014time} and \cite{yong2017linear} develop a discrete approximation to the value functions of time-inconsistent LQ problems in continuous time setting. \cite{dou2020time} extend this approximation to the Hilbert space using the semi-group method in partial differential equations. \cite{huang2021strong} and \cite{he2021equilibrium}  elaborate the feedback type equilibria by introducing the notions of ``strong equilibria'' and ``regular equilibria'' respectively. A different definition of equilibria is based on the notion of open loop control, which is proposed by \cite{hu2012time}. Moreover, \cite{hu2017time} derive a general necessary and sufficient condition for the open loop equilibrium for a time-inconsistent LQ problem and then proves the uniqueness of such type of equilibria in the one dimensional case\footnote{ Similar to \cite{basak2010dynamic} and \cite{bjomurzho2014}, the time inconsistency of the LQ problem discussed in \cite{hu2012time} and \cite{hu2017time} arises from the quadratic term of the expected state.}. \cite{yan2019open} obtain the existence and uniqueness of open loop equilibria for a mean-variance portfolio selection problem under stochastic volatility. Recently, the time-inconsistent LQ control models have been applied to other areas, such as mean-field games and differential games. The readers could be referred to \cite{bensoussan2013linear}, \cite{moon2020linear}, \cite{ni2017time}, \cite{lazrak2023time} and the reference theirin. Finally, it is worth noting that a non-existence result for general time-inconsistent control problems have been reported in literature. For a time-inconsistent binary control problem with non-exponential discounting, \cite{tan2021failure} find that an equilibrium may not exists, while the time-consistent counterpart admits a unique optimal solution. \par

The remainder  of the paper is organized as follows. Section \ref{Sec:ProblemSetting} introduces the formulation of the time-inconsistent LQ problem, in which we search for feedback intra-personal equilibria. Section \ref{Sec:EquilibriumPrinciple} derives the Bellman principle and HJB equation for the time-inconsistent LQ problem The precise definition of the solutions to the HJB equation is given.
Section \ref{SecUniqueness} characterizes the solutions to the HJB equation. Making use of the characterization of the value function of the LQ problem, we then obtain the main result: the equilibrium HJB equation admits at most one solution. Section \ref{Sec:Conclusion} concludes.

\section{Problem setting}\label{Sec:ProblemSetting}
For any $T>0$, we introduce the following notations.
\[
L^{p}\left((0,T);\mathbb{R}^{l\times k}\right) = \left\{f:(0,T)\to\mathbb{R}^{l\times k}\vert \int^T_0 \vert f_{ij}(t)\vert^p dt<\infty, 1\le i\le l, 1\le j\le k\right\}.
\]
\[
L^{\infty}\left((0,T);\mathbb{R}^{l\times k}\right) =\left \{f:(0,T)\to\mathbb{R}^{l\times k}\vert \esssup_{t\in[0,T]}| f_{ij}(t)|<\infty, 1\le i\le l, 1\le j\le k\right\}.
\]
\[C\left([0,T]^m; \mathbb{R}^{l\times k}\right)=\left\{f:[0,T]^m\to\mathbb{R}^{l\times k}\vert \text{$f$ is continuous}\right\}.\]
\begin{eqnarray*}
 C^{1}\left([0,T]^m; \mathbb{R}^{l\times k}\right)
 =\left\{f:[0,T]^m\to\mathbb{R}^{l\times k}\vert \text{$f$ is continuously differentiable}\right\}.
\end{eqnarray*}
\begin{eqnarray*}
 W^1\left((0,T);\mathbb{R}^{l\times k}\right) = \left\{f:(0,T)\to\mathbb{R}^{l\times k}\vert \text{ $f$ is weakly differentiable}\right\}.
\end{eqnarray*}
\begin{eqnarray*}
 WC^{1,1}\left((0,T)\times\mathbb{R}^{l\times k};\mathbb{R}\right) &=&\{f:(0,T)\times\mathbb{R}^{l\times k}\to\mathbb{R}\vert  f \mbox{ is weakly differentiable in }t \\
&&\mbox{and continuously differentiable in }x\}.
\end{eqnarray*}
\begin{eqnarray*}
 WC^{1,2}\left((0,T)\times\mathbb{R}^{l\times k};\mathbb{R}\right) &=&\{f:(0,T)\times\mathbb{R}^{l\times k}\to\mathbb{R}\vert  f\mbox{ is weakly differentiable in }t\\
  &&\mbox{and twice continuously differentiable in }x\}.
\end{eqnarray*}
For $f\in WC^{1,2}\left((0,T)\times\mathbb{R}^{l\times k};\mathbb{R}\right)$, $f_t$,  $\nabla_x f$ and $D_x^2 f$ denote the weak derivative in $t$, first order derivative and second order derivative in $x$ respectively.\par

For a real matrix-valued function $O(t)=\left(o_{ij}(t)\right)\in \mathbb{R}^{l\times k},$ $\forall t\in [0,T]^{m} (m=1,2)$, we introduce the following norms,
\begin{eqnarray*}
\left\{\begin{array}{ll}
\|O(t)\|=\max\limits_{1\leq i\leq l}\sum\limits_{j=1}^k|o_{ij}(t)|,\;\;\qquad
\|O\|_{L^p}=\max\limits_{1\leq i\leq l}\sum\limits_{j=1}^k\|o_{ij}\|_{L^p([0,T]^{m})}, (p=1,2)\\
\|O\|_{C}=\max\limits_{1\leq i\leq l}\sum\limits_{j=1}^k\|o_{ij}\|_{C([0,T]^{m})},\;\;\|O\|_{C^1}=\max\limits_{1\leq i\leq l}\sum\limits_{j=1}^k\left(\|o_{ij}\|_{C([0,T]^{m})}+\|Do_{ij}\|_{C([0,T]^{m})}\right).
\end{array}\right.
\end{eqnarray*}
For any initial pair $(t,x)\in[0,T)\times \mathbb{R}^n$, we consider the following controlled linear ordinary differential equation (LODE, for short)\footnote{Without any specification, any vector in this paper is a column vector.}
\begin{equation} \label{EDE}
\left\{
\begin{array}{ll}
\dot{y}(s)=A(s)y(s)+B(s)u(s)+b(s),& s\in(t,T], \\
y(t)=x,
\end{array}
\right.
\end{equation}
where the function  $u\in \mathcal{U}[0,T]  \equiv L^2 \left([0,T];\mathbb{R}^m\right)$ is the control process and $y$ is the state process valued in $\mathbb{R}^n.$

We suppose that the coefficients of the equation satisfy the following assumption
\begin{itemize}
\item[(H0)] $A\in L^{1}\left((0,T);\mathbb{R}^{n\times n}\right)$, $B\in L^{2}\left((0,T);\mathbb{R}^{n\times m}\right)$, $b\in L^1\left((0,T);\mathbb{R}^n\right)$.
\end{itemize}
It follows standard ODE theory that LODE (\ref{EDE}) has a unique solution $y(\cdot)\equiv  y^{u}_{t,x}(\cdot)$ in $C([t,T];\mathbb{R}^n)$, given by
\begin{align*}
y^{u}_{t,x} (s)=E(s,t)x+\int_t^s E(s,\tau)[B(\tau)u(\tau)+b(\tau)]d\tau,
\end{align*}
where
\[E(s,t) = \exp\left( \displaystyle\int ^s_tA(\tau)d\tau\right), \forall t,s\in[0,T].\]
At any time $t$ with the system state $y^{u}_{t,x}(t) = x$, the cost functional is given by
\begin{eqnarray} \label{CostFunctional}
J(t,x;u)&=&\int_{t}^{T}\bigg[\left\langle Q(t,s)y^{u}_{t,x}(s),y^{u}_{t,x}(s)\right\rangle+2\left\langle
S(t,s)y^{u}_{t,x}(s),u(s)\right\rangle\nonumber\\
&&+\left\langle M(t,s)u(s),u(s)\right\rangle
+2\left\langle q(t,s),y^{u}_{t,x}(s)\right\rangle+2\left\langle \rho(t,s),u(s)\right\rangle\bigg]ds\\
&&+\left\langle G(t)y^{u}_{t,x}(T),y^{u}_{t,x}(T)\right\rangle+2\left\langle g(t),y^{u}_{t,x}(T)\right\rangle,\nonumber
\end{eqnarray}
where the parameters satisfy the following assumptions throughout the paper:
\begin{itemize}
\item[(H1)]  $M\in C\left([0,T]\times[0,T]; \mathbb{R}^{m\times m}\right)$ is a positive  definite symmetric matrix-valued function.
\item[(H2)]  $Q\in C\left([0,T]\times[0,T];\mathbb{R}^{n\times n}\right) $ and $G\in C^1\left([0,T];\mathbb{R}^{n\times n}\right) $ are   positive semi-definite symmetric matrix-valued functions.
\item[(H3)]  $S\in C\left([0,T]\times[0,T]; \mathbb{R}^{m\times n}\right)$.
\item[(H4)] $q\in C\left([0,T]\times[0,T];\mathbb{R}^n\right)$, $\rho\in C\left([0,T]\times[0,T];\mathbb{R}^m\right)$, $g\in C^1\left([0,T];\mathbb{R}^n\right)$.
\item[(H5)] $Q$, $S$, $M$, $q$, $\rho$ are continuously differentiable with respect to the first variable and the corresponding derivatives are denoted by $Q_t$, $S_t$, $M_t$, $q_t$, $\rho_t$, respectively.
\end{itemize}
\begin{definition} A measurable mapping $\nu: [0,T]\times \mathbb{R}^n\to\mathbb{R}^m$ is an admissible feedback control law if the controlled evolution equation (\ref{EDE}) associated with $\nu$ admits a unique solution for any initial pair $(t,x)$. The set of all admissible feedback control laws is denoted by $\mathcal{A}$.
\end{definition}\par
For ease of exposition, we by abuse of notation denote the controlled system with initial pair $(t,x)$ and cost functional associated with the control law $\nu$ by $y^{\nu}$ and $J(t,x;\nu), \forall \nu\in\mathcal{A}$ respectively throughout the rest of the paper.
\begin{definition}[\cite{bjomur2009}]\label{DefEquilibrium}
The measurable mapping $\bar{u}:[0,T]\times \mathbb{R}^n\to \mathbb{R}^m$ is an equilibrium control law if $\bar{u}\in\mathcal{A}$ and
\begin{equation} \label{2.3}
 \liminf_{\varepsilon\searrow 0} \frac{J\left(t,x;u^{\varepsilon,v}\right)-J(t,x;\bar{u})}{\varepsilon}\geq 0,\;\;\forall (t,x,v)\in [0, T )\times \mathbb{R}^n\times\mathbb{R}^m,
 \end{equation}
 where
 \begin{equation} \label{2.4}
u^{\varepsilon,v}(s,y)=\left\{\begin{array}{ll}
v,&s\in [t,t+\varepsilon),y\in\mathbb{R}^n,\\
\bar{u}(s,y),&s\in [t+\varepsilon,T], y\in \mathbb{R}^n.
\end{array}
\right.
\end{equation}
\end{definition}
\begin{remark}
Definition \ref{DefEquilibrium} is established by seminal paper \cite{bjomur2009}, which states that for an equilibrium control law $\hat{u}$, any local perturbation will not make the performance functional better off. It is worth noting that the perturbation acts on the control law rather than the control process. This feature gives rise to a significant difference from the conventional maximum principle, in which the action is taken on control process, rather than the control law. In literature on time-inconsistent control, the two definitions would lead to vastly different equilibria. See \cite{he2022time} for more discussions.
\end{remark}
\section{The equilibrium Bellman principle and HJB equations}\label{Sec:EquilibriumPrinciple}
\subsection{The equilibrium Bellman principle}
In order to present the equilibrium HJB equation associated with equilibrium control law $\bar{u}$, we introduce the equilibrium value function $V$ and the equilibrium error function $R$ as follows:
\begin{eqnarray}\label{2.5}
V(t,x)=J(t,x;\bar{u}),\;\;\forall  (t,x)\in [0,T]\times \mathbb{R}^n
\end{eqnarray}
and
\begin{eqnarray}\label{2.6}
R(t,x)&=&\left\langle \dot{G}(t)y_{t,x}^{\bar{u}}(T)+2 \dot{g}(t),y_{t,x}^{\bar{u}}(T)\right\rangle+\int_t^T\bigg[\left\langle Q_t(t,s)y_{t,x}^{\bar{u}}(s)+2 q_t(t,s),y_{t,x}^{\bar{u}}(s)\right\rangle\nonumber\\
&&+\left\langle M_t(t,s)\bar{u}\left(s,y_{t,x}^{\bar{u}}(s)\right)+2 S_t(t,s)y_{t,x}^{\bar{u}}(s)+2 \rho_t(t,s),\bar{u}\left(s,y_{t,x}^{\bar{u}}(s)\right)\right\rangle\bigg]ds,
\end{eqnarray}
$\forall (t,x)\in [0,T]\times \mathbb{R}^n$.
\begin{theorem}\label{peng-theorem2.1} Suppose that Assumptions (H0)-(H5) hold and the equilibrium value function $V\in WC^{1,1}$ $\left([0,T]\times \mathbb{R}^n;\mathbb{R}\right)$,  then the equilibrium value function $V$ satisfies the following equilibrium Bellman principle
\begin{eqnarray}\label{2.7}
V(t,x)&=&\inf\limits_{u\in L^2((t,s);\mathbb{R}^m)}\bigg[\int_t^s\bigg(\left\langle Q(\tau,\tau)y_{t,x}^u(\tau)+2q(\tau,\tau),y_{t,x}^u(\tau)\right\rangle+2\left\langle S(\tau,\tau)y_{t,x}^u(\tau),u(\tau)\right\rangle\nonumber\\
&&+\langle M(\tau,\tau)u(\tau)+2\rho(\tau,\tau),u(\tau)\rangle-R\left(\tau,y_{t,x}^u(\tau)\right)\bigg)d\tau+V\left(s,y_{t,x}^u(s)\right)\bigg],
\end{eqnarray}
$\forall (t,x)\in [0,T]\times \mathbb{R}^n$ and $s\in[t,T]$.
\end{theorem}
\begin{proof}
Let us first consider the running cost $\displaystyle\int_{t}^{T}\left\langle Q(t,s)y^{u}_{t,x}(s),y^{u}_{t,x}(s)\right\rangle ds$.
\begin{eqnarray*}
&&\displaystyle\int_{t}^{T}\left\langle Q(t,s)y^{u}_{t,x}(s),y^{u}_{t,x}(s)\right\rangle ds\\
&=&\displaystyle\int_{t}^{s}\left\langle Q(t,\tau)y^{u}_{t,x}(\tau),y^{u}_{t,x}(\tau)\right\rangle d\tau+ \displaystyle\int_{s}^{T}\left\langle Q(t,\tau)y^{u}_{t,x}(\tau),y^{u}_{t,x}(\tau)\right\rangle d\tau\\
& =& \displaystyle\int_{t}^{s}\left\langle Q(\tau,\tau)y^{u}_{t,x}(\tau),y^{u}_{t,x}(\tau)\right\rangle d\tau - \displaystyle\int_{t}^{s}\displaystyle\int_{t}^{\tau}\left\langle Q_{\nu}(\nu,\tau)y^{u}_{t,x}(\tau),y^{u}_{t,x}(\tau)\right\rangle d\nu d\tau\\
&&+\displaystyle\int_{s}^{T}\left\langle Q(s,\tau)y^{u}_{t,x}(\tau),y^{u}_{t,x}(\tau)\right\rangle d\tau - \displaystyle\int_{s}^{T}\displaystyle\int_{t}^{s}\left\langle Q_{\nu}(\nu,\tau)y^{u}_{t,x}(\tau),y^{u}_{t,x}(\tau)\right\rangle d\nu d\tau.
\end{eqnarray*}
Changing the order of the double integration, we have that
\begin{eqnarray*}
&&\displaystyle\int_{t}^{T}\left\langle Q(t,s)y^{u}_{t,x}(s),y^{u}_{t,x}(s)\right\rangle ds\\
&=& \displaystyle\int_{t}^{s}\left\langle Q(\tau,\tau)y^{u}_{t,x}(\tau),y^{u}_{t,x}(\tau)\right\rangle d\tau - \displaystyle\int_{t}^{s}\displaystyle\int_{\nu}^{s}\left\langle Q_{\nu}(\nu,\tau)y^{u}_{t,x}(\tau),y^{u}_{t,x}(\tau)\right\rangle d\tau d\nu\\
&&+\displaystyle\int_{s}^{T}\left\langle Q(s,\tau)y^{u}_{t,x}(\tau),y^{u}_{t,x}(\tau)\right\rangle d\tau - \displaystyle\int_{t}^{s}\displaystyle\int_{s}^{T}\left\langle Q_{\nu}(\nu,\tau)y^{u}_{t,x}(\tau),y^{u}_{t,x}(\tau)\right\rangle d\tau d\nu \\
&=& \displaystyle\int_{t}^{s}\left\langle Q(\tau,\tau)y^{u}_{t,x}(\tau),y^{u}_{t,x}(\tau)\right\rangle d\tau -\displaystyle\int_{t}^{s}\displaystyle\int_{\nu}^{T}\left\langle Q_{\nu}(\nu,\tau)y^{u}_{t,x}(\tau),y^{u}_{t,x}(\tau)\right\rangle d\tau d\nu\\
&&+\displaystyle\int_{s}^{T}\left\langle Q(s,\tau)y^{u}_{t,x}(\tau),y^{u}_{t,x}(\tau)\right\rangle d\tau.
\end{eqnarray*}
Moreover, for the terminal pay-off $\left\langle G(t)y^{u}_{t,x}(T),y^{u}_{t,x}(T)\right\rangle$, we have that
\begin{align*}
\left\langle G(t)y^{u}_{t,x}(T),y^{u}_{t,x}(T)\right\rangle = \left\langle G(s)y^{u}_{t,x}(T),y^{u}_{t,x}(T)\right\rangle- \displaystyle\int^s_t\left\langle \dot{G}(\tau)y^{u}_{t,x}(T),y^{u}_{t,x}(T)\right\rangle d\tau.
\end{align*}
Let $\bar{u}$ denote an equilibrium control law. Using the similar calculations as the above on all the running costs and terminal pay-offs of $V$ and combining the representation of $V$ and $R$ , i.e., (\ref{2.5}) and (\ref{2.6}), we then have
\begin{eqnarray}\label{2.9}
V(t,x)&=&\int_{t}^{s}\bigg[\left\langle Q(\tau,\tau)y^{\bar{u}}_{t,x}(\tau),y^{\bar{u}}_{t,x}(\tau)\right\rangle+2\left\langle
S(\tau,\tau)y^{\bar{u}}_{t,x}(\tau),\bar{u}\left(\tau,y^{\bar{u}}_{t,x}(\tau)\right)\right\rangle\nonumber\\
&&+\left\langle M(\tau,\tau)\bar{u}\left(\tau,y^{\bar{u}}_{t,x}(\tau)\right),\bar{u}\left(\tau,y^{\bar{u}}_{t,x}(\tau)\right)\right\rangle+2\left\langle q(\tau,\tau),y^{\bar{u}}_{t,x}(\tau)\right\rangle\\
&&+2\left\langle \rho(\tau,\tau),\bar{u}\left(\tau,y^{\bar{u}}_{t,x}(\tau)\right)\right\rangle-R\left(\tau,y^{\bar{u}}_{t,x}(\tau)\right)\bigg]d\tau
+V\left(s,y^{\bar{u}}_{t,x}(s)\right).\nonumber
\end{eqnarray}

Next, we consider the perturbation control law $u^{\varepsilon,v}$ given by (\ref{2.4}). Solving the control system (\ref{EDE}) with $u^{\varepsilon,v}$, we have that the controlled systems has a  solution $ y_{t,x}^{u^{\varepsilon,v}}\in C([t,T];\mathbb{R}^n)$ given by
	\begin{equation} \label{2.11}
	y_{t,x}^{u^{\varepsilon,v}}(s)=\left\{\begin{array}{ll}
	y^{v}_{t,x}(s),&s\in [t,t+\varepsilon),\\
	y^{\bar{u}}_{t+\varepsilon,y^{v}_{t,x}(t+\varepsilon)}(s),&s\in [t+\varepsilon,T]
	\end{array}
	\right.
	\end{equation}
and
\begin{equation} \label{2.12}
\lim\limits_{\varepsilon \rightarrow0}y_{t,x}^{u^{\varepsilon,v}}=y_{t,x}^{\bar{u}} \mbox{ in } C\left([t,T];\mathbb{R}^n\right).
\end{equation}
Similar to (\ref{2.9}), applying change of the order of integrations to $J\left(t,x;u^{\varepsilon,v}\right)$, we then have
\begin{eqnarray*}
&&J\left(t,x;u^{\varepsilon,v}\right)\\
&=&\int_{t}^{t+\varepsilon}\bigg[\left\langle Q(t,s)y^{v}_{t,x}(s)+2q(t,s),y^{v}_{t,x}(s)\right\rangle+\left\langle M(t,s)v+2 S(t,s)y^{v}_{t,x}(s)+2\rho(t,s),v\right\rangle\bigg]ds\\
&&-\int_{t+\varepsilon}^{T}\int_{t}^{t+\varepsilon}\bigg[\left\langle Q_\tau(\tau,s)y^{\bar{u}}_{t+\varepsilon,y^{v}_{t,x}(t+\varepsilon)}(s)+2q_\tau(\tau,s),y^{\bar{u}}_{t+\varepsilon,y^{v}_{t,x}(t+\varepsilon)}(s)\right\rangle \\
&&+2\left\langle   S_\tau(\tau,s)y^{\bar{u}}_{t+\varepsilon,y^{v}_{t,x}(t+\varepsilon)}(s)+\rho_\tau(\tau,s),
\bar{u}\left(s,y^{\bar{u}}_{t+\varepsilon,y^{v}_{t,x}(t+\varepsilon)}(s)\right)\right\rangle \\
&&+\left \langle
 M_\tau(\tau,s)\bar{u}\left(s,y^{\bar{u}}_{t+\varepsilon,y^{v}_{t,x}(t+\varepsilon)}(s)\right),\bar{u}\left(s,y^{\bar{u}}_{t+\varepsilon,y^{v}_{t,x}(t+\varepsilon)}(s)\right)\right\rangle \bigg]d\tau ds\\
&&-\int_{t}^{t+\varepsilon}\left[\left\langle \dot{G}(\tau)y^{\bar{u}}_{t+\varepsilon,y^{v}_{t,x}(t+\varepsilon)}(T)+2\dot{g}(\tau),y^{\bar{u}}_{t+\varepsilon,y^{v}_{t,x}(t+\varepsilon)}(T)\right\rangle\right]d\tau+V\left(t+\varepsilon,y^{v}_{t,x}(t+\varepsilon)\right).
\end{eqnarray*}
It follows from (\ref{2.6}), (\ref{2.11}) and (\ref{2.12}) that
\begin{eqnarray}\label{variation}
&&\lim_{\varepsilon\searrow 0} \frac{J\left(t,x;u^{\varepsilon,v}\right)- J\left(t,x;\bar{u}\right)}{\varepsilon}\nonumber\\
&=&\left\langle Q(t,t)x+2q(t,t),x\right\rangle+2\left\langle S(t,t)x,v\right\rangle+\left\langle M(t,t)v+2\rho(t,t),v\right\rangle\nonumber\\
&&\quad-\int_{t}^{T}\bigg[\left\langle Q_t(t,s)y^{\bar{u}}_{t,x}(s)+2 q_t(t,s),y^{\bar{u}}_{t,x}(s)\right\rangle+2\left\langle   S_t(t,s)y^{\bar{u}}_{t,x}(s),\bar{u}\left(s,y^{\bar{u}}_{t,x}(s)\right)\right\rangle\nonumber\\
&&\quad+\left \langle M_t(t,s)\bar{u}\left(s,y^{\bar{u}}_{t,x}(s)\right),\bar{u}\left(s,y^{\bar{u}}_{t,x}(s)\right)\right\rangle +2\left\langle   \rho_t(t,s),\bar{u}\left(s,y^{\bar{u}}_{t,x}(s)\right)\right\rangle\bigg]ds\nonumber\\
&&\quad-\left\langle \dot{G}(t)y^{\bar{u}}_{t,x}(T)+2\dot{g}(t),y^{\bar{u}}_{t,x}(T)\right\rangle+V_t\left(t,x\right)+\langle \nabla_x V\left(t,x\right),A(t)x+B(t)v+b(t)\rangle\nonumber\\
&&=V_t\left(t,x\right)+\langle  \nabla_x V\left(t,x\right),A(t)x+B(t)v+b(t)\rangle+\langle Q(t,t)x,x\rangle+2\langle S(t,t)x,v\rangle\nonumber\\
&&\quad+\langle M(t,t)v,v\rangle+2\left\langle q(t,t),x\right\rangle+2\left\langle \rho(t,t),v\right\rangle-R(t,x).
\end{eqnarray}
Then the definition of equilibrium (\ref{2.3}) yields that
\begin{eqnarray*}\label{2.13}
&&V_t\left(t,x\right)+\langle \nabla_x V\left(t,x\right),A(t)x+B(t)v+b(t)\rangle+\langle Q(t,t)x,x\rangle+2\langle S(t,t)x,v\rangle\nonumber\\
&&+\langle M(t,t)v,v\rangle+2\left\langle q(t,t),x\right\rangle+2\left\langle \rho(t,t),v\right\rangle-R(t,x)\geq0,
\end{eqnarray*}
which implies that
\begin{eqnarray*}
&&\frac{d}{d\tau}V\left(\tau,y^{u}_{t,x}(\tau)\right)+\left\langle Q(\tau,\tau)y^{u}_{t,x}(\tau),y^{u}_{t,x}(\tau)\right\rangle+2\left\langle S(\tau,\tau)y^{u}_{t,x}(\tau),u(\tau)\right\rangle+\langle M(\tau,\tau)u(\tau),u(\tau)\rangle\nonumber\\
&&+2\left\langle q(\tau,\tau),y^{u}_{t,x}(\tau)\right\rangle+2\left\langle \rho(\tau,\tau),u(\tau)\right\rangle-R\left(\tau,y^{u}_{t,x}(\tau)\right)\geq0, \forall\tau\in [t,s],u\in L^2\left((t,s);\mathbb{R}^m\right).
\end{eqnarray*}
Integrating the above inequality from $t$ to $s$, we obtain
\begin{eqnarray*}
V(t,x)&\leq&\int_t^s\bigg[\left\langle Q(\tau,\tau)y^{u}_{t,x}(\tau),y^{u}_{t,x}(\tau)\right\rangle+2\left\langle S(\tau,\tau)y^{u}_{t,x}(\tau),u(\tau)\right\rangle+\langle M(\tau,\tau)u(\tau),u(\tau)\rangle\nonumber\\
&&+2\left\langle q(\tau,\tau),y^{u}_{t,x}(\tau)\right\rangle+2\left\langle \rho(\tau,\tau),u(\tau)\right\rangle-R\left(\tau,y^{u}_{t,x}(\tau)\right)\bigg]d\tau+V\left(s,y^{u}_{t,x}(s)\right).
\end{eqnarray*}
This and (\ref{2.9}) imply that (\ref{2.7}) holds and complete the proof.
\end{proof}
\subsection{The equilibrium HJB equation and verification}
Intuitively, Theorem~\ref{peng-theorem2.1} implies that the $V,R,\bar{u}$ satisfy the following equilibrium HJB equation
\begin{eqnarray}\label{2.8}
\left\{\begin{array}{ll}
V_t(t,x)+\left\langle \nabla_x V(t,x),  A(t)x\right\rangle+\inf\limits_{v\in \mathbb{R}^m}H(t,x,\nabla_x V(t,x),v)=0,\\
\qquad\qquad\qquad\qquad\qquad\qquad\qquad\qquad\mbox{   }  (t,x)\in[0,T)\times \mathbb{R}^n,\\
\bar{u}(t,x) =\argmin\limits_{v\in\mathbb{R}^m} H(t,x,\nabla_x V(t,x),v),(t,x)\in [0,T]\times \mathbb{R}^n,\\
V(T,x)=\langle G(T)x+2g(T), x\rangle,\qquad\quad x\in \mathbb{R}^n,
\end{array}\right.
\end{eqnarray}
where the Hamiltonian $H:[0,T]\times \mathbb{R}^n\times \mathbb{R}^n\times \mathbb{R}^m\rightarrow \mathbb{R}$ is defined by
\begin{eqnarray}\label{2.15}
H(t,x,p,v)&=&\left\langle p,B(t)v+b(t)\right\rangle+\left\langle Q(t, t)x,x\right\rangle+2\left\langle S(t, t)x,v \right\rangle+\left\langle M(t, t) v,v \right\rangle \nonumber\\
&&+2\left\langle q(t, t),x\right\rangle+2\left\langle \rho(t, t),v\right\rangle-R(t,x),(t,x,p,v)\in [0,T]\times \mathbb{R}^n\times \mathbb{R}^n\times \mathbb{R}^m.
\end{eqnarray}
In this subsection, we verify that the solutions of the equilibrium HJB equation (\ref{2.8}) indeed solve the equilibrium problem in Definition \ref{DefEquilibrium}.
\begin{theorem}\label{peng-theorem2.2}Suppose that Assumptions (H0)-(H5) hold. If the equilibrium HJB equation  (\ref{2.8}) admits a solution $V\in WC^{1,2}\left([0,T]\times \mathbb{R}^n;\mathbb{R}\right)$ and there exists $C>0$ such that 
\[\left\|D^2_xV(t,x)\right\|<C, \forall(t,x)\in [0,T]\times \mathbb{R}^n,\]
 then there exists an equilibrium control law given by
\begin{eqnarray}\label{2.16}
\bar{u}(t,x)=-M^{-1}(t,t)\left(\frac{1}{2}B^\top(t)\nabla_x V(t,x)+S(t,t)x+\rho(t,t)\right),
\end{eqnarray}
$\forall (t,x)\in [0,T]\times \mathbb{R}^n.$ Moreover,
\begin{align*}
V(t,x) = J(t,x;\bar{u}),\forall (t,x)\in [0,T]\times \mathbb{R}^n.
\end{align*}
\end{theorem}
\begin{proof}
Plugging the feedback control law $\bar{u}$ given by (\ref{2.16}) into the controlled systems (\ref{EDE}) and we have
\[
	\left\{
	\begin{array}{ll}
	\dot{\bar{y}}(s)=A(s)\bar{y}(s)+b(s)\\
\qquad\quad-B(s)M^{-1}(s,s)\left[\frac{1}{2}B^\top(s)\nabla_x V(s,\bar{y}(s))+S(s,s)\bar{y}(s)+\rho(s,s)\right], s\in(t,T], \\
	\bar{y}(t)=x.
	\end{array}
	\right.
\]
Given the regularity of $\nabla_x V$, the above Cauchy problem can be re rewritten as 
\[
	\left\{
	\begin{array}{ll}
	\dot{\bar{y}}(s)=\left[A(s)-B(s)M^{-1}(s,s)\left(\frac{1}{2}B^\top(s)\int_0^1D^2_x V(s,\theta\bar{y}(s))d\theta+S(s,s)\right)\right]\bar{y}(s)\\
\qquad\quad+b(s)+B(s)M^{-1}(s,s)\left[\frac{1}{2}B^\top(s)\nabla_x V(s,0)-\rho(s,s)\right],\quad s\in(t,T], \\
	\bar{y}(t)=x.
	\end{array}
	\right.
\]
We then have that the above ODE admits a unique solution $\bar{y}=y_{t,x}^{\bar{u}}\in C\left([t,T];\mathbb{R}^n\right)$, $\forall (t,x)\in [0,T)\times \mathbb{R}^n$.

For any given $(t,x)\in [0,T]\times \mathbb{R}^n$, it follows from (\ref{2.15}) that $H(t,x,\nabla_x V(t,x),\cdot)$ is a strictly convex function on $\mathbb{R}^m$ and $\bar{u}(t,x)$ is the unique minimum point of $H(t,x,\nabla_x V(t,x),\cdot)$. Following (\ref{2.6}) and (\ref{2.8}),  we then have
\begin{eqnarray*} \label{2.17}
\left\{\begin{array}{ll}
\frac{d}{ds}V(s,\bar{y}(s))+\left\langle M(s,s) \bar{u}(s,\bar{y}(s))+2\rho(s,s),\bar{u}(s,\bar{y}(s))\right\rangle +2\left\langle S(s, s)\bar{y}(s),\bar{u}(s,\bar{y}(s))\right\rangle\\
\quad +\left\langle Q(s, s)\bar{y}(s)+2q(s,s),\bar{y}(s)\right\rangle-\int_s^T\bigg[\left\langle M_s(s,\tau)\bar{u}\left(\tau,\bar{y}(\tau)\right)+2\rho_s(s,\tau),\bar{u}\left(\tau,\bar{y}(\tau)\right)\right\rangle\\
\quad+2\left\langle S_s(s,\tau)\bar{y}(\tau),\bar{u}\left(\tau,\bar{y}(\tau)\right)\right\rangle+\left\langle Q_s(s,\tau)\bar{y}(\tau)+2q_s(s,\tau),\bar{y}(\tau)\right\rangle\bigg]d\tau\\
\quad-\left\langle \dot{G}(s)\bar{y}(T)+2\dot{g}(s),\bar{y}(T)\right\rangle=0,\qquad s\in [t, T),\\
V(T,\bar{y}(T))=\langle G(T)\bar{y}(T)+2g(T),\bar{y}(T)\rangle.
\end{array}\right.
\end{eqnarray*}
Integrating the above equation from $t$ to $T$, we obtain
\begin{eqnarray*}
V(t,x)&=&\int_t^T\bigg[\left\langle M(s,s) \bar{u}(s,\bar{y}(s))+2\rho(s,s),\bar{u}(s,\bar{y}(s))\right\rangle+2\left\langle S(s, s)\bar{y}(s),\bar{u}(s,\bar{y}(s))\right\rangle\\
&&+\left\langle Q(s, s)\bar{y}(s)+2q(s,s),\bar{y}(s)\right\rangle\bigg]ds +\left\langle G(t)\bar{y}(T)+2g(t),\bar{y}(T)\right\rangle\\
&&-\int_t^T\int_s^T\bigg[\left\langle M_s(s,\tau)\bar{u}\left(\tau,\bar{u}(\tau)\right)+2\rho_s(s,\tau),\bar{u}\left(\tau,\bar{y}(\tau)\right)\right\rangle\\
&&+2\left\langle S_s(s,\tau)\bar{y}(\tau),\bar{u}\left(\tau,\bar{y}(\tau)\right)\right\rangle
+\left\langle Q_s(s,\tau)\bar{y}(\tau)+2q_s(s,\tau),\bar{y}(\tau)\right\rangle\bigg]d\tau ds\\
&=&\int_t^T\bigg[\left\langle M(s,s) \bar{u}(s,\bar{y}(s))+2\rho(s,s),\bar{u}(s,\bar{y}(s))\right\rangle+2\left\langle S(s, s)\bar{y}(s),\bar{u}(s,\bar{y}(s))\right\rangle\\
&&+\left\langle Q(s, s)\bar{y}(s)+2q(s,s),\bar{y}(s)\right\rangle\bigg]ds +\left\langle G(t)\bar{y}(T)+2g(t),\bar{y}(T)\right\rangle\\
&&-\int_t^T\int_t^\tau\bigg[\left\langle M_s(s,\tau)\bar{u}\left(\tau,\bar{y}(\tau)\right)+2\rho_s(s,\tau),\bar{u}\left(\tau,\bar{y}(\tau)\right)\right\rangle\\
&&+2\left\langle S_s(s,\tau)\bar{y}(\tau),\bar{u}\left(\tau,\bar{y}(\tau)\right)\right\rangle+\left\langle Q_s(s,\tau)\bar{y}(\tau)+2q_s(s,\tau),\bar{y}(\tau)\right\rangle\bigg] ds d\tau,
\end{eqnarray*}
which gives that
\begin{eqnarray*} \label{2.18}
V(t,x)&=&\int_t^T\bigg[\left\langle Q(t, s)\bar{y}(s)+2q(t,s),\bar{y}(s)\right\rangle+2\left\langle S(t, s)\bar{y}(s),\bar{u}(s,\bar{y}(s))\right\rangle\nonumber\\
&&+\left\langle M(t,s) \bar{u}(s,\bar{y}(s))+2\rho(t,s),\bar{u}(s,\bar{y}(s))\right\rangle\bigg]ds +\left\langle G(t)\bar{y}(T)+2g(t),\bar{y}(T)\right\rangle\nonumber\\
&=&J(t,x;\bar{u}),\forall (t,x)\in [0,T]\times \mathbb{R}^n.
\end{eqnarray*}

It suffices to verify $\bar{u}$ given by (\ref{2.16}) is an equilibrium control law. This is an immediate result of (\ref{variation}) and the HJB equation (\ref{2.8}), which completes the proof.
\end{proof}
\begin{remark} In the proof of Theorem \ref{peng-theorem2.2}, we can see that the regularity of $V$ ensures the existence and uniqueness of the solution to the controlled system (\ref{EDE}), and thus guarantees the admissibility of $\bar{u}$.
\end{remark}\par
We are now in the position to introduce a solution to the equilibrium HJB equation (\ref{2.8}) formally. To this end,
we introduce the following functions.
\begin{eqnarray*}\label{wp1}
\left\{\begin{array}{ll}
h(t,x,p)=M^{-1}(t,t)\left(\frac{1}{2}B^\top(t) p+S(t,t)x+\rho(t,t)\right),\\
\tilde{H}(t,x,p)=\left\langle \frac{1}{2}B^\top(t)p-S(t,t)x-\rho(t,t),h(t,x,p) \right\rangle+\left\langle Q(t, t)x+2q(t,t),x\right\rangle,\\
F(t,s,x,p)=\left\langle M_t(t,s)h(s,x,p)-2 S_t(t,s)x-2 \rho_t(t,s), h(s,x,p)\right\rangle\\
\qquad\qquad\quad\qquad+\left\langle Q_t(t,s)x+2 q_t(t,s),x\right\rangle,
\end{array}\right.
\end{eqnarray*}
$\forall (t,x,p)\in [0,T]\times \mathbb{R}^n\times \mathbb{R}^n$ and $s\in[t,T]$. \par

Replace $\bar{u}$ by the representation (\ref{2.16}) in the equilibrium HJB equation (\ref{2.8}), then the equilibrium HJB equation can be rewritten as
\begin{eqnarray}\label{wp2}
\left\{\begin{array}{ll}
V_t(t,x)+\left\langle \nabla_x V(t,x),  A(t)x-B(t)h(t,x, \nabla_x V(t,x))+b(t)\right\rangle\\
\qquad+\tilde{H}(t,x,\nabla_x V(t,x))-\int_t^TF\left(t,s,Y(s),\nabla_x V\left(s,Y(s)\right)\right)ds\\
\qquad-\left\langle \dot{G}(t)Y(T)+2 \dot{g}(t),Y(T)\right\rangle=0,\qquad\qquad\quad(t,x)\in [0,T]\times \mathbb{R}^n,\\
\dot{Y}(s)=A(s)Y(s)-B(s)h(s,Y(s), \nabla_x V(s,y(s)))+b(s),\quad  s\in (t,T],\\
V(T,x)=\langle G(T)x+2g(T), x\rangle, Y(t)=x, \qquad\qquad(t,x)\in [0,T]\times \mathbb{R}^n.
\end{array}\right.
\end{eqnarray}
We use the integral form of (\ref{wp2}) as the formal definition of the solution to the equilibrium HJB equation as follows.
\begin{definition}\label{peng-d2.2} The function $(Y,V)\in C([0,T];\mathbb{R}^n)\times (WC^{1,1}((0,T)\times \mathbb{R}^n;\mathbb{R})\cap C([0,T]\times \mathbb{R}^n;\mathbb{R}))$ is a solution to the equilibrium HJB equation (\ref{wp2})(i.e. (\ref{2.8})), if $(Y,V)$ satisfies the following integral equations
\begin{eqnarray*}
\left\{\begin{array}{ll}
V(t,x)=\left\langle G(t)Y(T)+2g(t),Y(T)\right\rangle+\displaystyle\int_t^T\bigg[\tilde{H}(\tau,Y(\tau),\nabla_x V(\tau,Y(\tau)))\\
\qquad\qquad-\int_\tau^TF\left(\tau,s,Y(s),\nabla_x V\left(s,Y(s)\right)\right)ds\bigg]d\tau,\\
Y(s)=E(s,t)x+\int_t^sE(s,\tau)\left[b(\tau)-B(\tau)h(\tau,Y(\tau), \nabla_x V(\tau,Y(\tau)))\right]d\tau,
\end{array}\right.
\end{eqnarray*}
$\forall (t,x)\in [0,T]\times \mathbb{R}^n$ and $t\leq s\leq T$.
\end{definition}
\section{Uniqueness of the solutions to the equilibrium HJB equation}\label{SecUniqueness}
In this section, we will discuss the uniqueness of solutions to the equilibrium HJB equation (\ref{2.8}). The uniqueness result is based on the uniqueness of solutions to a Riccati equation and the structure of the equilibrium value function.
\subsection{A Riccati equation}
In order to show the uniqueness result for the equilibrium HJB equation (\ref{2.8}), we introduce the  following equilibrium Riccati equation,
\begin{eqnarray} \label{3.1}
\left\{\begin{array}{ll}
 \dot{P}(t)+A^\top(t)P(t)+P(t)A(t)+ Q(t,t)-\mathbb{Q}(t)-\Gamma^\top(t)M(t,t)\Gamma(t)=0,\\
 \qquad\qquad\qquad\qquad\qquad\qquad\qquad\qquad\qquad\qquad\qquad\qquad\qquad t\in [0, T ),\\
 P(T)= G(T),
\end{array}\right.
\end{eqnarray}
where
 \begin{eqnarray}\label{3.2}
\left\{\begin{array}{ll}
\Gamma(t)=M^{-1}(t,t)\left(B^\top(t)P(t)+S(t,t)\right),\\
\mathbb{Q}(t)=\mathbb{E}^\top(T,t)\dot{G}(t)\mathbb{E}(T,t)+\int_t^T\mathbb{E}^T(s,t)\bigg[ Q_t(t,s)-\Gamma^\top(s)S_t(t,s)\\
\qquad\quad -S^\top_t(t,s)\Gamma(s)+\Gamma^\top(s)M_t(t,s)\Gamma(s)\bigg]\mathbb{E}(s,t)ds,
\end{array}
\right.
\end{eqnarray}
and
\begin{eqnarray*}\label{3.3}
\mathbb{E}(t,s)=E(t,s)-\int_s^tE(t,\tau)B(\tau)\Gamma(\tau)\mathbb{E}(\tau,s)d\tau,
\end{eqnarray*}
or equivalently,
\begin{align}\label{fundamentalSolution}
\mathbb{E}(t,s) = \exp\left(\int^t_s\left(A(\tau)-B(\tau)M^{-1}(\tau,\tau)\left(B^\top(\tau)P(\tau)+S(\tau,\tau)\right)\right)d\tau\right),
\end{align}
$\forall 0\leq s\leq t\leq T$.
\begin{definition}\label{SolutionRiccati}
$P\in C([0,T];\mathbb{R}^{n\times n})$  is a solution of the equilibrium Riccati equation (\ref{3.1}) if $P$ has the following integral form
\begin{eqnarray*}\label{3.4}
P(t)&=&\mathbb{E}^\top(T,t)G(T)\mathbb{E}(T,t)
+\int_{t}^T\mathbb{E}^\top(\tau,t)\bigg[Q(\tau,\tau)-S^\top(\tau,\tau)M^{-1}(\tau,\tau)S(\tau,\tau)\nonumber\\
&&-\mathbb{Q}(\tau)+P(\tau)B(\tau)M^{-1}(\tau,\tau)B^{\top}(\tau)P(\tau)\bigg]\mathbb{E}(\tau,t)d\tau,\quad t\in[0,T],
\end{eqnarray*}
 or equivalently
\begin{eqnarray}\label{3.5}
P(t)&=&E^\top(T,t)G(T)E(T,t)\nonumber\\
&+&\int_{t}^TE^\top(\tau,t)\bigg[Q(\tau,\tau)-\mathbb{Q}(\tau)-\Gamma^\top(\tau)M(\tau,\tau)\Gamma(\tau)\bigg]E(\tau,t)d\tau,\quad t\in[0,T].
\end{eqnarray}
\end{definition}
\begin{proposition}\label{pl1}
If that Assumptions (H0)-(H5) hold and $P_1, P_2\in C([0,T]$; $\mathbb{R}^{n\times n})$ are solutions to the Riccati equation (\ref{3.1}), then $P_1(t) = P_2(t)$ for all $ t\in[0,T]$.
\end{proposition}
\begin{proof}
Let $P_1$ and $P_2$ denote the two solutions to the equilibrium differential Riccati equation (\ref{3.1}). We introduce operator $\mathbb{G}$ as follows
\begin{eqnarray}\label{4.22}
\mathbb{G}(s;t,P_i)&=&\Psi_i^\top(T,s)\dot{G}(t)\Psi_i(T,s)+\int_{s}^T\Psi_i^\top(\tau,s)\Gamma_i^\top(\tau)\frac{\partial}{\partial t} M(t,\tau)\Gamma_i(\tau)\Psi_i(\tau,s)  d\tau
\nonumber\\
&&+\int_{s}^T\Psi_i^\top(\tau,s)\left[\frac{\partial}{\partial t} Q(t,\tau)-\frac{\partial}{\partial t}S^\top(t,\tau)\Gamma_i(\tau)-\Gamma_i^\top(\tau)\frac{\partial}{\partial t}S(t,\tau)\right] \Psi_i(\tau,s) d\tau,
\end{eqnarray}
$\forall s\in[0,T]$, where
\begin{align}\label{w0}
\Gamma_i(t)=M^{-1}(t,t)\left(B^\top(t)P_i(t)+S(t,t)\right), \forall t\in [0,T],
\end{align}
\begin{eqnarray}\label{w1}
\Psi_i(t,s)=E(t,s)-\int_s^tE(t,\tau)B(\tau)\Gamma_i(\tau)\Psi_i(\tau,s)d\tau,\quad \forall 0\leq s\leq t\leq T.
\end{eqnarray}
Thanks to (H0) and (\ref{w0}), we have
\begin{eqnarray}\label{w3}
\left\{\begin{array}{ll}
\|E(t,s)\|_{\mathbb{R}^{n\times n}}=\left\|\exp\left( \displaystyle\int ^s_tA(\tau)d\tau\right)\right\|_{\mathbb{R}^{n\times n}} \le \lambda,&0\leq t\leq s\leq T,\\
\int_t^s\|\Gamma_i(t)\|_{\mathbb{R}^{m\times n}}^2dt\leq\lambda_0^2\left(\|P_i\|_{C}+\|S\|_{C}\right)^2,&0\leq t\leq s\leq T,
\end{array}\right.
\end{eqnarray}
 where
\[\lambda = e^{\|A\|_{L^1}}, \lambda_0=\left\|M^{-1}\right\|_{C}\left(1+T+\|B\|_{L^2}\right).\]
This, (\ref{w0}) and (\ref{w1}) imply that
\begin{eqnarray*}
\|\Psi_i(t,s)\|_{\mathbb{R}^{n\times n}}\leq \lambda +\lambda\int_t^s \|B(\tau)\|_{\mathbb{R}^{n\times m}}\|\Gamma_i(t)\|_{\mathbb{R}^{m\times n}}\|\Psi_i(\tau,t)\|_{\mathbb{R}^{n\times n}}d\tau.
\end{eqnarray*}
Then it follows from H$\ddot{o}$lder's inequality and Gronwall's inequality that
\begin{eqnarray} \label{wp11}
\|\Psi_i(t,s)\|_{\mathbb{R}^{n\times n}}\leq \lambda\exp\left(\lambda\lambda_0\|B\|_{L^2}\left(\|P_i\|_{C}+\|S\|_{C}\right)\right),\forall t,s\in [0,T], i=1,2.
\end{eqnarray}
 Moreover, estimating $\mathbb{G}$ by (\ref{4.22}), we have that
 \begin{eqnarray}\label{wp22}
\|\mathbb{G}(s;t,P_i)\|_{\mathbb{R}^{n\times n}}
&\leq&\|G\|_{C^1}\|\Psi_i(T,s)\|^2_{\mathbb{R}^{n\times n}}\nonumber\\
&&+\tilde{\beta}\int_{s}^T\|\Psi_i(\tau,s)\|^2_{\mathbb{R}^{n\times n}}\left[ 1+\left\|\Gamma_i(\tau)\right\|_{\mathbb{R}^{m\times n}}\right]^2 d\tau\nonumber\\
&\leq&\lambda^2\tilde{\alpha}\left[ 1+T+\left\|\Gamma_i\right\|_{L^2}\right]^2\exp\left(2\lambda\lambda_0\|B\|_{L^2}\left(\|P_i\|_{C}+\|S\|_{C}\right)\right) \nonumber\\
&\leq&\alpha_0\left(1+\|P_i\|_{C}+\|S\|_{C}\right)^2\exp\left(2\lambda\lambda_0\|B\|_{L^2}\left(\|P_i\|_{C}+\|S\|_{C}\right)\right),
 \end{eqnarray}
$\forall s\in[0,T]$, where
\[\tilde{\beta}=\|Q\|_{C^1}+\|M\|_{C^1}+\|S\|_{C^1}\|,\tilde{\alpha}=\|G\|_{C^1}+  \tilde{\beta},\alpha_0=\lambda^2\left(1+T+\lambda_0\right)^2\tilde{\alpha}.\]
(\ref{wp11}) and (\ref{wp22}) tell that there exist constants $L(\|P_i\|_C)>0, i=1,2,$ such that
\begin{align}\label{4.22p}
\max\{\|\Psi_i(s,t)\|_{\mathbb{R}^{n\times n}},\|\mathbb{G}(s;t,P_i)\|_{\mathbb{R}^{n\times n}}\}\leq L(\|P_i\|_C),\quad 0\leq t\leq s\leq T.
\end{align}

For any $t\in[0,T]$, we define
\begin{eqnarray*}
\mathcal{L}(s)=\mathbb{G}(s;t,P_2)-\mathbb{G}(s;t,P_1) \quad \forall s\in[t,T].
\end{eqnarray*}
Then it follows from (\ref{4.22}) that $\mathcal{L}\in C\left([t,T];\mathbb{R}^n\right)\bigcap W^{1}\left((t,T);\mathbb{R}^n\right)$ satisfies the following Liyapulov equation
\begin{eqnarray*}
\left\{\begin{array}{ll}
\dot{\mathcal{L}}(s)+\left(A(s)-B(s)\Gamma_2(s)\right)^{\top}\mathcal{L}(s)+\mathcal{L}(s)\left(A(s)-B(s)\Gamma_2(s)\right)\\
\qquad+\left(P_1(s)-P_2(s)\right)B(s)M^{-1}(s,s)\left(B^\top(s)\mathbb{G}(s;t,P_1)-M_t(t,s) \Gamma_2(s)+2S_t(t,s)\right)\\
\qquad+\left(\mathbb{G}(s;t,P_1)B(s)-\Gamma^{\top}_1(s)M_t(t,s)\right)M^{-1}(s,s)B^\top(s)\left(P_1(s)-P_2(s)\right)=0,\\
\mathcal{L}(T)=0,
\end{array}\right.
\end{eqnarray*}
$\forall s\in[t,T)$.  Solving the above Liyapulov equation, we have that
\begin{eqnarray*}
\mathcal{L}(s)&=&\mathbb{G}(s;t,P_2)-\mathbb{G}(s;t,P_1)\nonumber\\
&=&\int_{s}^{T}\Psi_2^{\top}(\tau,s)\bigg[\left(\mathbb{G}(\tau;t,P_1)B(\tau)-\Gamma^{\top}_1(\tau)M_t(t,\tau)\right)M^{-1}(\tau,\tau)B^{\top}(\tau)\left(P_1(\tau)-P_2(\tau)\right)\\
&&+\left(P_1(\tau)-P_2(\tau)\right)B(\tau)M^{-1}(\tau,\tau)\bigg(B^{\top}(\tau)\mathbb{G}(\tau;t,P_1)\\
&&-M_t(t,\tau) \Gamma_2(\tau)+2S_t(t,\tau)\bigg)\bigg]\Psi_2(\tau,s)d\tau.\nonumber
\end{eqnarray*}
Plugging the estimate (\ref{4.22p}) into the above integral, we have that
\begin{eqnarray}\label{w4}
&&\|\mathbb{G}(t;t,P_2)-\mathbb{G}(t;t,P_1)\|_{\mathbb{R}^{n\times n}}\nonumber\\
&\leq&2L^2(\|P_2\|_C)\left\|M^{-1}\right\|_C\int_{t}^{T}\bigg[L(\|P_1\|_C)\|B(\tau)\|_{\mathbb{R}^{n\times m}}+\|M\|_{C^1}\|\Gamma_1(\tau)\|_{\mathbb{R}^{m\times n}}\nonumber\\
&&+\|M\|_{C^1}\| \Gamma_2(\tau)\|_{\mathbb{R}^{m\times n}}+\|S\|_{C^1}\bigg]\|B(\tau)\|_{\mathbb{R}^{n\times m}}\left\|P_1(\tau)-P_2(\tau)\right\|_{\mathbb{R}^{n\times n}}d\tau\nonumber\\
&\leq&\beta_0\int_{t}^{T}\bigg[1+\|B(\tau)\|_{\mathbb{R}^{n\times m}}+\|\Gamma_1(\tau)\|_{\mathbb{R}^{m\times n}}\nonumber\\
&&+\| \Gamma_2(\tau)\|_{\mathbb{R}^{m\times n}}\bigg]\|B(\tau)\|_{\mathbb{R}^{n\times m}}\left\|P_1(\tau)-P_2(\tau)\right\|_{\mathbb{R}^{n\times n}}d\tau,
\end{eqnarray}
where
\[\beta_0=2L^2(\|P_2\|_C)\left\|M^{-1}\right\|_C\left(L(\|P_1\|_C)+\|\|M\|_{C^1}+\|S\|_{C^1}\right).\]

On the other hand, following (\ref{3.5}), we have that
\begin{eqnarray}\label{4.24}
P_1(t)-P_2(t)=\int_{t}^TE^\top(\tau,t)\big[\mathbb{G}(\tau;\tau,P_2)-\mathbb{G}(\tau;\tau,P_1)\qquad\qquad\qquad\qquad\qquad\qquad\nonumber\\
+\left(P_2(\tau)-P_1(\tau)\right)B(\tau)\Gamma_2(\tau)+\Gamma_1^\top(\tau)B^{\top}(\tau)\left(P_2(\tau)-P_1(\tau)\right)\big]E(\tau,t)d\tau,
\end{eqnarray}
$\forall t\in[0,T]$. 
Plugging (\ref{w3}) and (\ref{w4}) into (\ref{4.24}), we obtain
\begin{eqnarray*}
\|P_1(t)-P_2(t)\|_{\mathbb{R}^{n\times n}}
&\leq&\gamma_0\int_{t}^T\bigg[1+\|B(\tau)\|_{\mathbb{R}^{n\times m}}+\|\Gamma_1(\tau)\|_{\mathbb{R}^{m\times n}}\\
&&+\| \Gamma_2(\tau)\|_{\mathbb{R}^{m\times n}}\bigg]\|B(\tau)\|_{\mathbb{R}^{n\times m}}\left\|P_1(\tau)-P_2(\tau)\right\|_{\mathbb{R}^{n\times n}}d\tau,
\end{eqnarray*}
where $\gamma_0=2(1+T\beta_0)\lambda^2$.  Together with (\ref{w0}), we have
\begin{eqnarray*}
\|P_1(t)-P_2(t)\|_{\mathbb{R}^{n\times n}}
&\leq&\gamma_0\left(1+\left\|M^{-1}\right\|_C\right)(\|P_1\|_{C}+\|P_2\|_C\\
&&+2\| S\|_{C})\int_{t}^T\left(1+\|B(\tau)\|_{\mathbb{R}^{n\times m}}\right)^2\left\|P_1(\tau)-P_2(\tau)\right\|_{\mathbb{R}^{n\times n}}d\tau.
\end{eqnarray*}
Then it follows from Gronwall's inequality that $P_2=P_1$ in $C\left([0,T];\mathbb{R}^{n\times n}\right)$, which completes the proof.
\end{proof}
\subsection{Two auxiliary functions}
We will see later that the equilibrium value function $V$ can be represented by $P$ and two auxiliary functions $\varphi$ and $\psi$, which satisfy the following two equations respectively,
\begin{eqnarray} \label{3.22}
\left\{\begin{array}{ll}
\dot{\varphi}(t)+\left(A(t)-B(t)\Gamma(t)\right)^\top \varphi(t)-\mathbb{S}(t) + P(t)b(t)\\
\qquad+q(t,t)- \Gamma^\top(t)\rho(t,t)=0,&t\in [0, T ),\\
\varphi(T)=g(T),
\end{array}\right.
\end{eqnarray}
\begin{eqnarray} \label{3.23}
\left\{\begin{array}{ll}
\dot{\psi}(t)+2\langle  \varphi(t),-B(t)\Upsilon(t)+b(t)\rangle-\omega(t)\\
\qquad+\langle M(t,t)\Upsilon(t)-2\rho(t,t),\Upsilon(t)\rangle=0,& t\in [0, T ),\\
\psi(T)=0,
\end{array}\right.
\end{eqnarray}
where,
\begin{eqnarray}\label{3.24}
\left\{\begin{array}{ll}
\Upsilon(t)=M^{-1}(t,t)\left(B^\top(t) \varphi(t)+\rho(t,t)\right),\\
\tilde{b}(s,t)=\int_t^s \mathbb{E}(s,\tau)(b(\tau)-B(\tau)\Upsilon(\tau))d\tau,\\
\mathbb{S}(t)=\int_t^T\mathbb{E}^\top(s,t)\bigg[\left(Q_t(t,s)-\Gamma^\top(s) S_t(t,s)-S^\top_t(t,s)\Gamma(s)\right)\tilde{b}(s,t)\\
\quad+\Gamma^\top(s)M_t(t,s)\left(\Upsilon(s)+\Gamma(s)\tilde{b}(s,t)\right)-\Gamma^\top(s)\rho_t(t,s)+ q_t(t,s)\\
\quad -S^\top_t(t,s)\Upsilon(s)\bigg]ds+\mathbb{E}^\top(T,t)\dot{g}(t)+\mathbb{E}^\top(T,t)\dot{G}(t)\tilde{b}(T,t),\\
\omega(t)=\left\langle \dot{G}(t)\tilde{b}(T,t)+2 \dot{g}(t),\tilde{b}(T,t)\right\rangle+\int_t^T\left\langle Q_t(t,s)\tilde{b}(s,t),\tilde{b}(s,t)\right\rangle ds\\
\quad+2\int_t^T\left\langle  q_t(t,s)-\Gamma^\top(s) S_t(t,s)\tilde{b}(s,t)-S^\top_t(t,s)\Upsilon(s),\tilde{b}(s,t)\right\rangle ds\\
\quad+\int_t^T\bigg\langle M_t(t,s)[\Gamma(s)\tilde{b}(s,t)+\Upsilon(s)]-2\rho_t(t,s),\Gamma(s)\tilde{b}(s,t)+\Upsilon(s)\bigg\rangle ds,
\end{array}
\right.
\end{eqnarray}
$0\le t\leq s\le T$. \par
Note that both (\ref{3.22}) and (\ref{3.23}) are LODEs, then it is easy to obtain the following result.
\begin{proposition}\label{wpl2} Suppose that Assumptions (H0)-(H5) hold. For any $P\in C\left([0,T];\mathbb{R}^{n\times n}\right)$, LODEs (\ref{3.22}) and (\ref{3.23}) admit unique solutions in $C([0,T])$.
\end{proposition}

\subsection{Characterization of the solution to the equilibrium HJB equation}

\begin{proposition}\label{pl2} Suppose that Assumptions (H0)-(H5) hold and $(Y,V)\in C([0$, $ T]; \mathbb{R}^n)\times (WC^{1,2}((0,T)\times \mathbb{R}^n;\mathbb{R})\cap C([0,T]\times \mathbb{R}^n;\mathbb{R}))$ is a solution to the equilibrium HJB equation (\ref{2.8}). Let $P(s) := \int^1_0\frac{1}{2} D^2_x V(s,\theta Y(s))d\theta.$ Then
\begin{itemize}
\item[(i)] $P$ is the solution to the equilibrium Riccati equation  (\ref{3.1}).
\item[(ii)] $P(s)$ is independent of the initial state $(t,x)$ of $Y(s), s\ge t$.
\item[(iii)] $\nabla_x V(t,x) = 2P(t)x + \varphi(t), \forall (t,x)\in [0,T]\times\mathbb{R}^n,$ where $\varphi$ is the solution to LODE (\ref{3.22}).
\end{itemize}
\end{proposition}
\begin{proof}
$(i)$
Let $y_{t,x}^{\bar{u}}\in C([t,T];\mathbb{R}^n)$ denote a solution of the controlled system (\ref{EDE}) with the equilibrium control $\bar{u}$. Then, it follows from (\ref{2.16}) that
\begin{eqnarray}\label{4.3}
\bar{u}\left(s,y_{t,x}^{\bar{u}}(s)\right)=-\tilde{\Gamma}(s)y_{t,x}^{\bar{u}}(s)-\tilde{\Upsilon}(s),s\in[t,T].
\end{eqnarray}
where
\begin{eqnarray}\label{4.2}
\left\{\begin{array}{ll}
\tilde{\Gamma}(s)= M^{-1}(s,s)\left(\frac{1}{2}B^\top(s)\int_0^1D^2_x V\left(s,\theta y_{t,x}^{\bar{u}}(s)\right)d\theta +S(s,s)\right),\\
\tilde{\Upsilon}(s)=M^{-1}(s,s)\left(\frac{1}{2}B^\top(s)\nabla_x V(s,0)+\rho(s,s)\right),
\end{array}\right.
\end{eqnarray}
for all $s\in[t,T]$.

Let
\begin{eqnarray}\label{4.4}
\mathbf{E}(\tau,s)=\exp\left(\int_s^{\tau}(A(\nu)-B(\mu)\tilde{\Gamma}(\nu))d\nu\right)\quad 0\leq t\leq s\leq \tau\leq T.
\end{eqnarray}
Hence, solving LODE (\ref{EDE}), we obtain that
\begin{eqnarray}\label{4.5}
y^{\bar{u}}_{t,x}(s)=\mathbf{E}(s,t)x+f(s,t),s\in[t,T],
\end{eqnarray}
where
\begin{eqnarray*}\label{4.6}
f(s,t)= \int_t^s\mathbf{E}(s,\tau)(b(\tau)-B(\tau)\tilde{\Upsilon}(\tau))d\tau.
\end{eqnarray*}
Plug (\ref{4.5}) into the definition of the equilibrium error function $R$ (see (\ref{2.6})), we then have
\begin{eqnarray*}
R(t,x)&=&\left\langle \dot{G}(t)(\mathbf{E}(T,t)x+f(T,t))+2 \dot{g}(t),\mathbf{E}(T,t)x+f(T,t)\right\rangle\\
&&+\int_t^T\bigg[\left\langle Q_t(t,s)(\mathbf{E}(s,t)x+f(s,t))+2 q_t(t,s),\mathbf{E}(s,t)x+f(s,t)\right\rangle\qquad\\
&&+\langle M_t(t,s)\bar{u}\left(s,\mathbf{E}(s,t)x+f(s,t)\right)+2 S_t(t,s)(\mathbf{E}(s,t)x+f(s,t))\\
&&+2 \rho_t(t,s),\bar{u}\left(s,\mathbf{E}(s,t)x+f(s,t)\right)\rangle\bigg]ds,\nonumber
\end{eqnarray*}
which yields
\begin{eqnarray}\label{4.7}
R(t,x)=\langle \tilde{\mathbb{Q}}(t)x,x\rangle+2\langle \tilde{\mathbb{S}}(t),x\rangle+\tilde{\omega}(t),
\end{eqnarray}
where
\begin{eqnarray}\label{4.8}
\left\{\begin{array}{ll}
\tilde{\mathbb{Q}}(t)=\mathbf{E}^\top(T,t)\dot{G}(t)\mathbf{E}(T,t)+\int_t^T\mathbf{E}^\top(s,t)\bigg[ Q_t(t,s)-\tilde{\Gamma}^\top(s)S_t(t,s)\\
\quad -S^\top_t(t,s)\tilde{\Gamma}(s)+\tilde{\Gamma}^\top(s)M_t(t,s)\tilde{\Gamma}(s)\bigg]\mathbf{E}(s,t)ds,\\
\tilde{\mathbb{S}}(t)=\mathbf{E}^\top(T,t)\dot{g}(t)+\mathbf{E}^\top(T,t)\dot{G}(t)f(T,t)\\
\quad+\int_t^T\mathbf{E}^\top(s,t)\bigg[Q_t(t,s)f(s,t)-\left(\tilde{\Gamma}^\top(s) S_t(t,s)+S^\top_t(t,s)\tilde{\Gamma}(s)\right)f(s,t)\\
\quad-S^\top_t(t,s)\tilde{\Upsilon}(s)+\tilde{\Gamma}^\top(s)M_t(t,s)\left(\tilde{\Upsilon}(s)+\tilde{\Gamma}(s)f(s,t)\right)\\
\quad+ q_t(t,s)-\tilde{\Gamma}^\top(s)\rho_t(t,s)\bigg]ds,\\
\tilde{\omega}(t)=\left\langle \dot{G}(t)f(T,t)+2 \dot{g}(t),f(T,t)\right\rangle+\int_t^T\left\langle Q_t(t,s)f(s,t),f(s,t)\right\rangle ds\\
\quad+2\int_t^T\left\langle  q_t(t,s)-\tilde{\Gamma}^\top(s) S_t(t,s)f(s,t)-S^\top_t(t,s)\tilde{\Upsilon}(s),f(s,t)\right\rangle ds\\
\quad+\int_t^T\bigg\langle M_t(t,s)[\tilde{\Gamma}(s)f(s,t)+\tilde{\Upsilon}(s)]-2\rho_t(t,s),\tilde{\Gamma}(s)f(s,t)+\tilde{\Upsilon}(s)\bigg\rangle ds.
\end{array}
\right.
\end{eqnarray}
Plug (\ref{4.3}) and (\ref{4.7}) into (\ref{2.8}), we have
\begin{eqnarray*}
&&V_t\left(s,y^{\bar u}_{t,x}(s)\right)+\left\langle \nabla_x V\left(s,y^{\bar u}_{t,x}(s)\right),  A(s)y^{\bar u}_{t,x}(s)\right\rangle+H\left(s,y^{\bar u}_{t,x}(s),\nabla_x V\left(s,y^{\bar u}_{t,x}(s)\right),\bar{u}\left(s,y^{\bar u}_{t,x}(s)\right)\right)\\
&=&\left\langle \nabla_x V\left(s,y^{\bar u}_{t,x}(s)\right),b(s)-B(s)\left(\tilde{\Gamma}(s)y^{\bar u}_{t,x}(s)+\tilde{\Upsilon}(s)\right)\right\rangle+\left\langle Q(s,s)y^{\bar u}_{t,x}(s),y^{\bar u}_{t,x}(s)\right\rangle-\tilde{\omega}(s)\\
&&+\left\langle M(s,s)\left(\tilde{\Gamma}(s)y^{\bar u}_{t,x}(s)+\tilde{\Upsilon}(s)\right),\tilde{\Gamma}(s)y^{\bar u}_{t,x}(s)+\tilde{\Upsilon}(s)\right\rangle+2\left\langle q(s,s),y^{\bar u}_{t,x}(s)\right\rangle -2\left\langle \tilde{\mathbb{S}}(s),y^{\bar u}_{t,x}(s)\right\rangle\\
&&-2\left\langle S(s,s)y^{\bar u}_{t,x}(s),\tilde{\Gamma}(s)y^{\bar u}_{t,x}(s)+\tilde{\Upsilon}(s) \right\rangle-2\left\langle \rho(s,s),\tilde{\Gamma}(s)y^{\bar u}_{t,x}(s)+\tilde{\Upsilon}(s)\right\rangle-\left\langle \tilde{\mathbb{Q}}(s)y^{\bar u}_{t,x}(s),y^{\bar u}_{t,x}(s)\right\rangle\\
&=&0\\
&\leq&V_t(s,\nu)+\left\langle \nabla_x V(s,\nu),  A(s)\nu\right\rangle+\left\langle \nabla_x V(s,\nu),b(s)-B(s)\left(\tilde{\Gamma}(s)\nu+\tilde{\Upsilon}(s)\right)\right\rangle+\left\langle Q(s,s)\nu,\nu\right\rangle\\
&&-2\left\langle S(s,s)\nu,\tilde{\Gamma}(s)\nu+\tilde{\Upsilon}(s) \right\rangle+\left\langle M(s,s)\left(\tilde{\Gamma}(s)\nu+\tilde{\Upsilon}(s)\right),\tilde{\Gamma}(s)\nu+\tilde{\Upsilon}(s)\right\rangle+2\left\langle q(s,s),\nu\right\rangle \nonumber\\
&&-2\left\langle \rho(s,s),\tilde{\Gamma}(s)\nu+\tilde{\Upsilon}(s)\right\rangle-\langle \tilde{\mathbb{Q}}(s)\nu,\nu\rangle-2\langle \tilde{\mathbb{S}}(s),\nu\rangle-\tilde{\omega}(s),\forall \nu\in \mathbb{R}^n.
\end{eqnarray*}
Using the optimality and taking the first order derivative on the above function (after the inequality sign of the above inequality) with respect to $\nu$, we then have
\begin{eqnarray*}
0&=&\frac{\partial}{\partial \nu}\bigg[V_t(s,\nu)+\left\langle \nabla_x V(s,\nu),  A(s)\nu\right\rangle+\left\langle \nabla_x V(s,\nu),b(s)-B(s)\left(\tilde{\Gamma}(s)\nu+\tilde{\Upsilon}(s)\right)\right\rangle+\left\langle Q(s,s)\nu,\nu\right\rangle\\
&&-2\left\langle S(s,s)\nu,\tilde{\Gamma}(s)\nu+\tilde{\Upsilon}(s) \right\rangle+\left\langle M(s,s)\left(\tilde{\Gamma}(s)\nu+\tilde{\Upsilon}(s)\right),\tilde{\Gamma}(s)\nu+\tilde{\Upsilon}(s)\right\rangle+2\left\langle q(s,s),\nu\right\rangle \nonumber\\
&&-2\left\langle \rho(s,s),\tilde{\Gamma}(s)\nu+\tilde{\Upsilon}(s)\right\rangle-\langle \tilde{\mathbb{Q}}(s)\nu,\nu\rangle-2\langle \tilde{\mathbb{S}}(s),\nu\rangle-\tilde{\omega}(s)\bigg]|_{\nu=y^{\bar u}_{t,x}(s)}\\
&=&\frac{d}{ds}V_x\left(s,y^{\bar u}_{t,x}(s)\right)+(A(s)-B(s)\tilde{\Gamma}(s))^\top\nabla_x V\left(s,y^{\bar u}_{t,x}(s)\right)\\
&&+2 \left[Q(s,s)- S^\top(s,s)\tilde{\Gamma}(s)-\tilde{\Gamma}^\top(s)S(s,s)+\tilde{\Gamma}^\top(s)M(s,s)\tilde{\Gamma}(s)-\tilde{\mathbb{Q}}(s)\right]y^{\bar u}_{t,x}(s)\\
&&+2 \tilde{\Gamma}^\top(s)M(s,s)\tilde{\Upsilon}(s)-2 S^\top(s,s)\tilde{\Upsilon}(s)+2 q(s,s)-2\tilde{\Gamma}^\top(s)\rho(s,s)-2 \tilde{\mathbb{S}}(s).
\end{eqnarray*}
Together with (\ref{2.5}) and (\ref{4.5}), then we have
\begin{eqnarray*}
\left\{\begin{array}{ll}
\frac{d}{ds}\nabla_x V\left(s,y_{t,x}^{\bar{u}}(s)\right)+(A(s)-B(s)\tilde{\Gamma}(s))^\top \nabla_x V\left(s,y_{t,x}^{\bar{u}}(s)\right)+2 q(s,s)\\
+2\left[Q(s,s)+\tilde{\Gamma}^\top(s)M(s,s)\tilde{\Gamma}(s)-\tilde{\Gamma}^\top(s)S(s,s)-S^\top(s,s)\tilde{\Gamma}(s)  -\tilde{\mathbb{Q}}(s) \right]y_{t,x}^{\bar{u}}(s) \\
-2S^\top(s,s)\tilde{\Upsilon}(s)-2\tilde{\Gamma}^\top(s) \rho(s,s)+2 \tilde{\Gamma}^\top(s)M(s,s)\tilde{\Upsilon}(s)-2\tilde{\mathbb{S}}(s)=0,\quad\forall s\in[t,T),\\
\nabla_x V(T,y_{t,x}^{\bar{u}}(T))=2 G(T)\mathbf{E}(T,t)x+2 G(T)f(T,t)+2g(T),
\end{array}
\right.
\end{eqnarray*}
Thus, the above  ODE with respect to  $\nabla _x V\left(\cdot,y_{t,x}^{\bar{u}}(\cdot)\right)\in C([t,T];\mathbb{R}^n)$ yields
\begin{eqnarray*}
\nabla_x V\left(s,y_{t,x}^{\bar{u}}(s)\right)
&=&2\int_s^T\mathbf{E}^\top(\tau,s)\bigg[Q(\tau,\tau)-\tilde{\mathbb{Q}}(\tau)+\tilde{\Gamma}^\top(\tau)M(\tau,\tau)\tilde{\Gamma}(\tau)-\tilde{\Gamma}^\top(\tau)S(\tau,\tau)\\
&&-S^\top(\tau,\tau)\tilde{\Gamma}(\tau)  \bigg]y_{t,x}^{\bar{u}}(\tau)d\tau+2\int_s^T\mathbf{E}^\top(\tau,s)\bigg[q(\tau,\tau)-S^\top(\tau,\tau)\tilde{\Upsilon}(\tau)\\
&&-\tilde{\Gamma}^\top(\tau) \rho(\tau,\tau)+ \tilde{\Gamma}^\top(\tau)M(\tau,\tau)\tilde{\Upsilon}(\tau)-\tilde{\mathbb{S}}(\tau)\bigg]d\tau\\
&&+2\mathbf{E}^\top(T,s)\left( G(T)\mathbf{E}(T,t)x+ G(T)f(T,t)+g(T)\right).
\end{eqnarray*}
Plug (\ref{4.5}) into the above equation,  we  then have
\begin{eqnarray}\label{New1}
\nabla_x V\left(s,y_{t,x}^{\bar{u}}(s)\right)
&=&2\mathbf{E}^\top(T,s)\left( G(T)\mathbf{E}(T,s)y_{t,x}^{\bar{u}}(s)+ G(T)f(T,s)+g(T)\right)\nonumber\\
&&+2\int_s^T\mathbf{E}^\top(\tau,s)\bigg[Q(\tau,\tau)-\tilde{\mathbb{Q}}(\tau)+\tilde{\Gamma}^\top(\tau)M(\tau,\tau)\tilde{\Gamma}(\tau)-\tilde{\Gamma}^\top(\tau)S(\tau,\tau)\nonumber\\
&&-S^\top(\tau,\tau)\tilde{\Gamma}(\tau)\bigg]\left(\mathbf{E}(\tau,s)y_{t,x}^{\bar{u}}(s)+ f(\tau,s)\right)d\tau+2\int_s^T\mathbf{E}^\top(\tau,s)\bigg[q(\tau,\tau)\\
&&+\tilde{\Gamma}^\top(\tau)M(\tau,\tau)\tilde{\Upsilon}(\tau)-S^\top(\tau,\tau)\tilde{\Upsilon}(\tau)-\tilde{\Gamma}^\top(\tau) \rho(\tau,\tau)-\tilde{\mathbb{S}}(\tau)\bigg]d\tau\nonumber,
\end{eqnarray}
$\forall s\in[t,T]$.\par

We define the  operator $P_{t,x}:[0,T]\longrightarrow \mathbb{R}^{n\times n}$ as follows
\begin{eqnarray} \label{4.15}
P_{t,x}(s)&=&\mathbf{E}^\top(T,s)G(T)\mathbf{E}(T,s)+\int_s^T\mathbf{E}^\top(\tau,s)\bigg[Q(\tau,\tau)-\tilde{\Gamma}^\top(\tau)S(\tau,\tau)\nonumber\\
&&+\tilde{\Gamma}^\top(\tau)M(\tau,\tau)\tilde{\Gamma}(\tau)-S^\top(\tau,\tau)\tilde{\Gamma}(\tau) -\tilde{\mathbb{Q}}(\tau)\bigg]\mathbf{E}(\tau,s)d\tau.
\end{eqnarray}
Then it is to see from (\ref{New1}) and (\ref{4.15}) that
\begin{align}\label{New2}
\nabla_x V\left(s,y_{t,x}^{\bar{u}}(s)\right) = 2P_{t,x}(s) y_{t,x}^{\bar{u}}(s)+\nabla_x V\left(s,0\right),\;\;s\in[t,T].
\end{align}
On the other hand,
\begin{align}\label{New3}
\nabla_x V\left(s,y_{t,x}^{\bar{u}}(s)\right) &=\int_0^1\frac{d}{d\theta}\nabla_x V\left(s,\theta y_{t,x}^{\bar{u}}(s)\right)d\theta+\nabla_x V\left(s,0\right)\nonumber\\
&=\int_0^1 D^2_x V\left(s,\theta y_{t,x}^{\bar{u}}(s)\right)d\theta y_{t,x}^{\bar{u}}(s)+\nabla_x V\left(s,0\right).
\end{align}
Compare (\ref{New2}) and (\ref{New3}), we then obtain
\begin{align}\label{New4}
2P_{t,x}(s) = \int_0^1 D^2_x V\left(s,\theta y_{t,x}^{\bar{u}}(s)\right)d\theta.
\end{align}

It suffices to show that $P_{t,x}$ solves the equilibrium Riccati equation (\ref{3.1}).
Plug (\ref{New4}) into (\ref{4.2}), then we have
\begin{eqnarray}\label{4.18}
\left\{\begin{array}{ll}
\tilde{\Gamma}(s)= M^{-1}(s,s)\left(B^\top(s)P_{t,x}(s) +S(s,s)\right),&s\in[t,T],\\
\tilde{\Upsilon}(s)=M^{-1}(s,s)\left(B^\top(s)\nabla_x V(s,0)+\rho(s,s)\right),&s\in[t,T].
\end{array}\right.
\end{eqnarray}
Plug the representation of $\tilde{\Gamma}$ in (\ref{4.18}) into (\ref{4.4}), we then have
\begin{eqnarray}\label{4.20}
\qquad\mathbf{E}(\tau,s)=\exp\left(\int_s^\tau \left(A(\nu)-B(\nu)M^{-1}(\nu,\nu)\left(B^\top(\nu)P_{t,x}(\nu) +S(\nu,\nu)\right)\right)d\nu\right),
\end{eqnarray}
$\forall 0\leq t\leq s\leq\tau\leq T$.\par 
It is to see from Assumptions (H0)-(H5) and (\ref{4.15}) that $P_{t,x}(s)$ is symmetric. Plug (\ref{4.18}) into (\ref{4.15}), we then have
\begin{eqnarray*} \label{4.19}
P_{t,x}(s)&=&\mathbf{E}^\top(T,s)G(T)\mathbf{E}(T,s)+\int_s^T\mathbf{E}^\top(\tau,s)\bigg[Q(\tau,\tau)-S^T(\tau,\tau)M^{-1}(\tau,\tau)S(\tau,\tau)\nonumber\\
&&+P_{t,x}(\tau)B(\tau)M^{-1}(\tau,\tau)B^\top(\tau)P_{t,x}(\tau)  -\tilde{\mathbb{Q}}(\tau)\bigg]\mathbf{E}(\tau,s)d\tau.
\end{eqnarray*}
 Similarly, plug (\ref{4.20}) and (\ref{4.4}) into $\tilde{\mathbb{Q}}$ (see (\ref{4.8})), we then obtain that $P_{t,x}$ is a symmetric solution to the Riccati equation (\ref{3.1}) in $C\left([0,T];\mathbb{R}^{n\times n}\right)$.

$(ii)$
Note that the coefficients of the Riccati equation (\ref{3.1}) are independent of $(t,x)$, hence it follows from Proposition \ref{pl1} that $P_{t,x}(s)$, i.e.,
\[P(s)=\int_0^1 \frac{1}{2}D^2_x V\left(s,\theta y_{t,x}^{\bar{u}}(s)\right)d\theta,\]
 is independent of $(t,x).$

$(iii)$  For any $(t,x)\in [0,T]\times \mathbb{R}^n$, it follows from (\ref{New3}) that
\begin{eqnarray*}
\nabla_x V\left(s, y^{\bar u}_{t,x}(s)\right)&=&\int_0^1 D_x^2 V\left(s, \theta y^{\bar u}_{t,x}(s)\right)d\theta y^{\bar u}_{t,x}(s)+ \nabla_x V\left(s, 0\right)\\
&=&2P(s) y^{\bar u}_{t,x}(s)+ \nabla_x V\left(s, 0\right),\forall s\in[t,T].
\end{eqnarray*}
Let $s=t$, then
\[\nabla_x V\left(t, x\right)=2P(t) x +\nabla_x V\left(t, 0\right), \forall(t,x)\in [0,T]\times \mathbb{R}^n.\]
It follows from (\ref{New1}) that
\begin{eqnarray*}
\nabla_x V\left(t,x\right)
&=&2\mathbf{E}^\top(T,t)\left( G(T)\mathbf{E}(T,t)x+ G(T)f(T,t)+g(T)\right)\nonumber\\
&&+2\int_t^T\mathbf{E}^\top(\tau,t)\bigg[Q(\tau,\tau)-\tilde{\mathbb{Q}}(\tau)+\tilde{\Gamma}^\top(\tau)M(\tau,\tau)\tilde{\Gamma}(\tau)-\tilde{\Gamma}^\top(\tau)S(\tau,\tau)\nonumber\\
&&-S^\top(\tau,\tau)\tilde{\Gamma}(\tau)\bigg]\left(\mathbf{E}(\tau,t)x+ f(\tau,t)\right)d\tau+2\int_t^T\mathbf{E}^\top(\tau,t)\bigg[q(\tau,\tau)\nonumber\\
&&+\tilde{\Gamma}^\top(\tau)M(\tau,\tau)\tilde{\Upsilon}(\tau)-S^\top(\tau,\tau)\tilde{\Upsilon}(\tau)-\tilde{\Gamma}^\top(\tau) \rho(\tau,\tau)-\tilde{\mathbb{S}}(\tau)\bigg]d\tau,t\in[0,T].\nonumber
\end{eqnarray*}
Then it is easy to verify that $\nabla_x V\left(t, 0\right)$ solves LODE (\ref{3.22}). This completes the proof.
\end{proof}

In the proof of Proposition \ref{pl2}, we introduce functions $\tilde{\mathbb{Q}}$, $\tilde{\mathbb{S}}$, $\tilde{\omega}$ in (\ref{4.8}). The following lemma shows that these functions are the same as functions $\mathbb{Q}$ in (\ref{3.2}), $\mathbb{S}$, $\omega$ in (\ref{3.24}) respectively. We will use this lemma in the proof of Theorem \ref{peng-theorem4.1}.
\begin{lemma}\label{Lem1}Suppose that Assumptions (H0)-(H5) hold and $(Y,V)\in C([0,T]$; $\mathbb{R}^n)\times (WC^{1,2}((0,T)\times \mathbb{R}^n;\mathbb{R})\cap C([0,T]\times \mathbb{R}^n;\mathbb{R}))$ is a solution to the equilibrium HJB equation (\ref{2.8}). Then
$\tilde{\mathbb{Q}}(t) = \mathbb{Q}(t), \tilde{\mathbb{S}}(t) = \mathbb{S}(t), \tilde{\omega}(t) = \omega(t),$ $\forall t\in [0,T].$
\end{lemma}
\begin{proof}
It follows from Proposition \ref{pl2} that
\[P(s) = \int^1_0\frac{1}{2} D^2_x V(s,\theta Y(s))d\theta\]
is the unique solution to the equilibrium Riccati equation. Then it is to see from (\ref{3.2}) and (\ref{4.2}) that
$\tilde{\Gamma}(t) = \Gamma(t), \forall t\in[0,T]$.  Thus, comparing (\ref{fundamentalSolution}) and (\ref{4.4}), we have that
\begin{align}\label{EandE}
\mathbb{E}(\tau,s) = \mathbf{E}(\tau,s), \quad \forall 0\le s\le \tau \le T.
\end{align}
Note that  $\nabla_x V\left(t, 0\right)$ solves LODE (\ref{3.22}), then Proposition \ref{wpl2}, (\ref{3.24}) and (\ref{4.2}) yield that
\begin{align}\label{UandU}
\tilde{\Upsilon}(t)=\Upsilon(t), \quad  t\in[0,T],
\end{align}
which, together with (\ref{EandE}), gives\footnote{See (\ref{4.6}) and (\ref{3.24}) for the definitions of $f$ and $\tilde{b}$ respectively.}
\begin{align}\label{FandB}
f(s,t) = \tilde{b}(s,t), \quad 0\le t\le s\le T.
\end{align}
Keeping (\ref{EandE}), (\ref{UandU}) and (\ref{FandB}) in mind and comparing (\ref{3.24}) and (\ref{4.8}), we have that $\tilde{\omega}(t)=\omega(t)$ and $\tilde{\mathbb{S}}(t) = \mathbb{S}(t), \forall t\in [0,T]$. Similarly, comparing (\ref{3.2}) and (\ref{4.8}), we have that $\tilde{\mathbb{Q}}(t) = \mathbb{Q}(t),\forall t\in[0,T].$ This completes the proof.
\end{proof}
\begin{theorem}\label{peng-theorem4.1}
Suppose that Assumptions (H0)-(H5) hold and $(Y,V)\in C([0,T]$; $\mathbb{R}^n)\times (WC^{1,2}((0,T)\times \mathbb{R}^n;\mathbb{R})\cap C([0,T]\times \mathbb{R}^n;\mathbb{R}))$ is a solution to the equilibrium HJB equation (\ref{2.8}). Then
\begin{eqnarray}\label{3.25}
V(t,x)=\langle P(t)x,x\rangle+2\langle\varphi(t),x\rangle+\psi(t), \quad (t,x)\in [0,T]\times \mathbb{R}^n,
\end{eqnarray}
\begin{align}\label{3.025}
Y(s) = \mathbb{E}(s,t)x + \int^s_t\mathbb{E}(s,\tau)\left[b(\tau)-B(\tau) M^{-1}(\tau,\tau)\left(B^\top(\tau)\varphi(\tau)+\rho(\tau,\tau)\right)\right]d\tau,
\end{align}
$\forall s\in[t,T]$, where $P$ is the solution to the equilibrium Riccati equation (\ref{3.1}), $\varphi$ and $\psi$ are the solutions to the (\ref{3.22}) and (\ref{3.23}) respectively, and $\mathbb{E}(s,t)$ is defined by (\ref{fundamentalSolution}).
\end{theorem}
\begin{proof}
It follows from Proposition \ref{pl2} that
\begin{eqnarray}\label{pw0}
V(t,x)=\langle P(t)x,x\rangle+2\langle \varphi(t),x\rangle+V(t,0),\quad \forall(t,x)\in [0,T]\times \mathbb{R}^n,
\end{eqnarray}
It suffices to show that $V(t,0)$ is the solution to LODE (\ref{3.23}). First, Lemma \ref{Lem1} and (\ref{4.7}) yield that
\begin{eqnarray}\label{4.21}
R(t,x)=\langle \mathbb{Q}(t)x,x\rangle+2\langle \mathbb{S}(t),x\rangle+\omega(t),\quad  (t,x)\in[0,T]\times \mathbb{R}^n.
\end{eqnarray}

Second, we define $V(t,0)$ by $\tilde{\psi}(t)$ and plug (\ref{pw0}) and (\ref{4.21}) into the HJB equation (\ref{2.8}), we then have
\begin{eqnarray*}
\left\{\begin{array}{ll}
0=\langle \dot{P}(t)x,x\rangle+2\langle \dot{\varphi}(t),x\rangle+\dot{\tilde{\psi}}(t)+2\left\langle P(t)x+\varphi(t),  A(t)x+B(t)\bar{u}(t,x)+b(t)\right\rangle\\
\qquad+\left\langle Q(t, t)x,x\right\rangle+2\left\langle S(t, t)x,\bar{u}(t,x)\right\rangle+\left\langle M(t, t)\bar{u}(t,x),\bar{u}(t,x) \right\rangle +2\left\langle q(t, t),x\right\rangle\\
\qquad+2\left\langle \rho(t, t),\bar{u}(t,x)\right\rangle-\langle \mathbb{Q}(t)x,x\rangle-2\langle \mathbb{S}(t),x\rangle-\omega(t),\qquad(t,x)\in [0,T]\times \mathbb{R}^n,\\
\langle P(T)x,x\rangle+2\langle \varphi(T),x\rangle+\tilde{\psi}(T)=\langle G(T)x+2g(T), x\rangle, x\in \mathbb{R}^n.
\end{array}\right.
\end{eqnarray*}
Let $x=0$, the above equation becomes
\begin{eqnarray*}
\left\{\begin{array}{ll}
0=\dot{\tilde{\psi}}(t)+2\left\langle \varphi(t),  B(t)\bar{u}(t,0)+b(t)\right\rangle+\left\langle M(t, t)\bar{u}(t,0),\bar{u}(t,0) \right\rangle \\
\qquad+2\left\langle \rho(t, t),\bar{u}(t,0)\right\rangle-\omega(t),&t\in [0,T),\\
\tilde{\psi}(T)=0.
\end{array}\right.
\end{eqnarray*}
By (\ref{4.3}) and (\ref{UandU}), we have
\begin{eqnarray*}
\bar{u}\left(t,0\right)=-\Upsilon(t),\;\; t\in[0,T].
\end{eqnarray*}
Moreover, we obtain
\begin{eqnarray*}
\left\{\begin{array}{ll}
0=\dot{\tilde{\psi}}(t)+2\left\langle \varphi(t), b(t) -B(t)\Upsilon(t)\right\rangle+\left\langle M(t, t)\Upsilon(t)-2 \rho(t, t),\Upsilon(t)\right\rangle-\omega(t),t\in [0,T),\\
\tilde{\psi}(T)=0.
\end{array}\right.
\end{eqnarray*}
Thus, $\tilde{\psi}$ solves LODE (\ref{3.23}). As LODE (\ref{3.23}) admits a unique solution in $C([0,T];\mathbb{R})$, determined by $P$ and $\varphi$, then Propositions \ref{pl1} and \ref{wpl2} yield that $\tilde{\psi}(t)=\psi(t)$ for all $t\in [0,T]$.

Finally, solving the ODE with respect to $Y$ in (\ref{wp2}), we obtain (\ref{3.025}). This completes the proof.
\end{proof}
\subsection{Uniqueness}
We are now in the position to show the uniqueness result for the equilibrium HJB equation  (\ref{2.8}).
\begin{theorem}\label{peng-corollary4.1}
Suppose that Assumptions (H0)-(H5) hold.
If $(Y_i,V_i)\in C([0,T]$; $\mathbb{R}^n)\times (WC^{1,2}((0,T)\times \mathbb{R}^n;\mathbb{R})\cap C([0,T]\times \mathbb{R}^n;\mathbb{R}))$ $(i=1,2)$ is a solution to the equilibrium HJB equation (\ref{2.8}). Then
\[
V_1(t,x) = V_2(t,x)\mbox{ and } Y_1(s) = Y_2(s),\quad  (t,x)\in [0,T]\times \mathbb{R}^n, s\in [t,T].
\]
\end{theorem}
\begin{proof}
Thanks to Theorem~\ref{peng-theorem4.1}, the uniqueness of the solutions to the equilibrium HJB equation (\ref{2.8}) boils down to the uniqueness of the solutions to the Riccati equation (\ref{3.1}) and LODEs (\ref{3.22}) and (\ref{3.23}), which are given by Propositions~\ref{pl1} and \ref{wpl2} respectively. This completes the proof.
\end{proof}

\section{Concluding remarks}\label{Sec:Conclusion}
The establishment of uniqueness of solutions to HJB equations holds pivotal significance within the realm of partial differential equation theory, particularly in its application to the study of control problems. We have obtained the uniqueness result for a class of HJB equations arising from general time-inconsistent deterministic LQ control problems. It is our hope that this result contributes to the broader understanding of general time-inconsistent control problems. Furthermore, we anticipate that the confirmed uniqueness, particularly in the context of our delineated class of HJB equations, will provide a foundational basis to validate, at least to some degree, the definition of equilibria across a spectrum of feedback controls. This, in turn, enhances the theoretical underpinning for the study and analysis of time-inconsistent control problems.

\section*{Acknowledgments}
We are grateful for comments from Xun Yu Zhou and Jiong Min Yong.

\bibliographystyle{ecta} 
\bibliography{SiamFinal}

\begin{thebibliography}{32}
\newcommand{\enquote}[1]{``#1''}
\expandafter\ifx\csname natexlab\endcsname\relax\def\natexlab#1{#1}\fi

\bibitem[\protect\citeauthoryear{Basak and Chabakauri}{Basak and
  Chabakauri}{2010}]{basak2010dynamic}
\textsc{Basak, S. and G.~Chabakauri} (2010): \enquote{Dynamic mean-variance
  asset allocation,} \emph{Review of financial Studies}, 23, 2970--3016.

\bibitem[\protect\citeauthoryear{Bensoussan, Sung, and Yam}{Bensoussan
  et~al.}{2013}]{bensoussan2013linear}
\textsc{Bensoussan, A., K.~Sung, and S.~C.~P. Yam} (2013):
  \enquote{Linear--quadratic time-inconsistent mean field games,} \emph{Dynamic
  Games and Applications}, 3, 537--552.

\bibitem[\protect\citeauthoryear{Bj\"{o}rk and Murgoci}{Bj\"{o}rk and
  Murgoci}{2009}]{bjomur2009}
\textsc{Bj\"{o}rk, T. and A.~Murgoci} (2009): \enquote{A general theory of
  Markovian time inconsistent stochastic control problems,} \emph{Working
  Paper}.

\bibitem[\protect\citeauthoryear{Bj\"{o}rk, Murgoci, and Zhou}{Bj\"{o}rk
  et~al.}{2014}]{bjomurzho2014}
\textsc{Bj\"{o}rk, T., A.~Murgoci, and X.~Zhou} (2014): \enquote{Mean-variance
  Portfolio Optimization with State Dependent Risk Aversion,}
  \emph{Mathematical Finance}, 24, 1--24.

\bibitem[\protect\citeauthoryear{Cai, Chen, Peng, and Wei}{Cai
  et~al.}{2022}]{cai2022time}
\textsc{Cai, H., D.~Chen, Y.~Peng, and W.~Wei} (2022): \enquote{On the
  time-inconsistent deterministic linear-quadratic control,} \emph{SIAM Journal
  on Control and Optimization}, 60, 968--991.

\bibitem[\protect\citeauthoryear{Dou and L\"{u}}{Dou and
  L\"{u}}{2020}]{dou2020time}
\textsc{Dou, F. and Q.~L\"{u}} (2020): \enquote{Time-inconsistent linear
  quadratic optimal control problems for stochastic evolution equations,}
  \emph{SIAM Journal on Control and Optimization}, 58, 485--509.

\bibitem[\protect\citeauthoryear{Ebert and Prelec}{Ebert and
  Prelec}{2007}]{ebepre2007}
\textsc{Ebert, J. and D.~Prelec} (2007): \enquote{The Fragility of Time:
  Time-Insensitivity and Valuation of the Near and Far Future,}
  \emph{Management Science}, 53, 1423--1438.

\bibitem[\protect\citeauthoryear{Ebert, Wei, and Zhou}{Ebert
  et~al.}{2020}]{ebert2020weighted}
\textsc{Ebert, S., W.~Wei, and X.~Y. Zhou} (2020): \enquote{Weighted
  discounting-on group diversity, time-inconsistency, and consequences for
  investment,} \emph{Journal of Economic Theory}, 189, 1--39.

\bibitem[\protect\citeauthoryear{Ekeland and Lazrak}{Ekeland and
  Lazrak}{2006}]{ekelaz2006}
\textsc{Ekeland, I. and A.~Lazrak} (2006): \enquote{Being Serious About
  Non-Commitment: Subgame Perfect Equilibrium in Continuous Time,} Working
  Paper.

\bibitem[\protect\citeauthoryear{Ekeland and Pirvu}{Ekeland and
  Pirvu}{2008}]{ekepir2008}
\textsc{Ekeland, I. and T.~Pirvu} (2008): \enquote{Investment and consumption
  without commitment,} \emph{Mathematics and Financial Economics}, 2, 57--86.

\bibitem[\protect\citeauthoryear{Harris and Laibson}{Harris and
  Laibson}{2013}]{harlai2013}
\textsc{Harris, C. and D.~Laibson} (2013): \enquote{Instantaneous
  Gratification,} \emph{Quarterly Journal of Economics}, 128, 205--248.

\bibitem[\protect\citeauthoryear{He and Jiang}{He and
  Jiang}{2021}]{he2021equilibrium}
\textsc{He, X.~D. and Z.~L. Jiang} (2021): \enquote{On the equilibrium
  strategies for time-inconsistent problems in continuous time,} \emph{SIAM
  Journal on Control and Optimization}, 59, 3860--3886.

\bibitem[\protect\citeauthoryear{He and Zhou}{He and Zhou}{2022}]{he2022time}
\textsc{He, X.~D. and X.~Y. Zhou} (2022): \enquote{Who are I: Time
  inconsistency and intrapersonal conflict and reconciliation,} in
  \emph{Stochastic Analysis, Filtering, and Stochastic Optimization: A
  Commemorative Volume to Honor Mark HA Davis's Contributions}, Springer,
  177--208.

\bibitem[\protect\citeauthoryear{Hu, Jin, and Zhou}{Hu
  et~al.}{2012}]{hu2012time}
\textsc{Hu, Y., H.~Jin, and X.~Y. Zhou} (2012): \enquote{Time-inconsistent
  stochastic linear--quadratic control,} \emph{SIAM journal on Control and
  Optimization}, 50, 1548--1572.

\bibitem[\protect\citeauthoryear{Hu, Jin, and Zhou}{Hu
  et~al.}{2017}]{hu2017time}
---\hspace{-.1pt}---\hspace{-.1pt}--- (2017): \enquote{Time-inconsistent
  stochastic linear-quadratic control: characterization and uniqueness of
  equilibrium,} \emph{SIAM Journal on Control and Optimization}, 55,
  1261--1279.

\bibitem[\protect\citeauthoryear{Huang and Zhou}{Huang and
  Zhou}{2021}]{huang2021strong}
\textsc{Huang, Y.-J. and Z.~Zhou} (2021): \enquote{Strong and weak equilibria
  for time-inconsistent stochastic control in continuous time,}
  \emph{Mathematics of Operations Research}, 46, 428--451.

\bibitem[\protect\citeauthoryear{Krussell and Smith}{Krussell and
  Smith}{2003}]{krusmi2003}
\textsc{Krussell, P. and A.~Smith} (2003): \enquote{Consumption-Savings
  Decision with Quasi-Geometric Discounting,} \emph{Econometrica}, 71,
  365--375.

\bibitem[\protect\citeauthoryear{Laibson}{Laibson}{1997}]{lai1997}
\textsc{Laibson, D.} (1997): \enquote{Golden Eggs and Hyperbolic Discounting,}
  \emph{Quarterly Journal of Economics}, 112, 443--378.

\bibitem[\protect\citeauthoryear{Lazrak, Wang, and Yong}{Lazrak
  et~al.}{2023}]{lazrak2023time}
\textsc{Lazrak, A., H.~Wang, and J.~Yong} (2023): \enquote{Time-Inconsistency
  in Linear Quadratic Stochastic Differential Games,} \emph{arXiv preprint
  arXiv:2304.11577}.

\bibitem[\protect\citeauthoryear{Luttmer and Mariotti}{Luttmer and
  Mariotti}{2003}]{lutmar2003}
\textsc{Luttmer, E. and T.~Mariotti} (2003): \enquote{Subjective Discounting in
  an Exchange Economy,} \emph{Journal of Politic Economy}, 11, 959--989.

\bibitem[\protect\citeauthoryear{Moon and Yang}{Moon and
  Yang}{2020}]{moon2020linear}
\textsc{Moon, J. and H.~J. Yang} (2020): \enquote{Linear-quadratic
  time-inconsistent mean-field type Stackelberg differential games:
  Time-consistent open-loop solutions,} \emph{IEEE Transactions on Automatic
  Control}, 66, 375--382.

\bibitem[\protect\citeauthoryear{Ni, Zhang, and Krstic}{Ni
  et~al.}{2017}]{ni2017time}
\textsc{Ni, Y.-H., J.-F. Zhang, and M.~Krstic} (2017):
  \enquote{Time-inconsistent mean-field stochastic LQ problem: Open-loop
  time-consistent control,} \emph{IEEE Transactions on Automatic Control}, 63,
  2771--2786.

\bibitem[\protect\citeauthoryear{O'Donoghue and Rabin}{O'Donoghue and
  Rabin}{2001}]{odorab1999}
\textsc{O'Donoghue, T. and M.~Rabin} (2001): \enquote{Choice and
  Procrastination,} \emph{Quarterly Journal of Economics}, 116, 112--160.

\bibitem[\protect\citeauthoryear{Phelps and Pollak}{Phelps and
  Pollak}{1968}]{phepol1968}
\textsc{Phelps, E. and R.~Pollak} (1968): \enquote{On Second-Best National
  Saving and Game-Equilibrium Growth,} \emph{Review of Economic Studies}, 35,
  185--199.

\bibitem[\protect\citeauthoryear{Strotz}{Strotz}{1955}]{str1955}
\textsc{Strotz, R.} (1955): \enquote{Myopia and Inconsistency in Dynamic
  Utility Maximization,} \emph{Review of Economic Studies}, 23, 165--180.

\bibitem[\protect\citeauthoryear{Tan, Wei, and Zhou}{Tan
  et~al.}{2021}]{tan2021failure}
\textsc{Tan, K.~S., W.~Wei, and X.~Y. Zhou} (2021): \enquote{Failure of smooth
  pasting principle and nonexistence of equilibrium stopping rules under
  time-inconsistency,} \emph{SIAM journal on control and optimization}, 59,
  4136--4154.

\bibitem[\protect\citeauthoryear{Thaler}{Thaler}{1981}]{tha1981}
\textsc{Thaler, R.} (1981): \enquote{Some Empirical Evidence on Dynamic
  Inconsistency,} \emph{Economics Letters}, 8, 201--207.

\bibitem[\protect\citeauthoryear{Weitzman}{Weitzman}{2001}]{weitzman2001gamma}
\textsc{Weitzman, M.~L.} (2001): \enquote{Gamma discounting,} \emph{American
  Economic Review}, 91, 260--271.

\bibitem[\protect\citeauthoryear{Yan and Wong}{Yan and
  Wong}{2019}]{yan2019open}
\textsc{Yan, T. and H.~Y. Wong} (2019): \enquote{Open-loop equilibrium strategy
  for mean--variance portfolio problem under stochastic volatility,}
  \emph{Automatica}, 107, 211--223.

\bibitem[\protect\citeauthoryear{Yong}{Yong}{2012}]{yong2012deterministic}
\textsc{Yong, J.} (2012): \enquote{A deterministic linear quadratic
  time-inconsistent optimal control problem,} \emph{arXiv preprint
  arXiv:1204.1856}.

\bibitem[\protect\citeauthoryear{Yong}{Yong}{2014}]{yong2014time}
---\hspace{-.1pt}---\hspace{-.1pt}--- (2014): \enquote{Time-inconsistent
  optimal control problems,} in \emph{Proceedings of the international congress
  of mathematicians}, vol.~4, 947--969.

\bibitem[\protect\citeauthoryear{Yong}{Yong}{2017}]{yong2017linear}
---\hspace{-.1pt}---\hspace{-.1pt}--- (2017): \enquote{Linear-quadratic optimal
  control problems for mean-field stochastic differential
  equations-time-consistent solutions,} \emph{Transactions of the American
  Mathematical Society}, 369, 5467--5523.

\end{thebibliography}
\end{document}